\documentclass[a4,10.5pt,leqno]{article}
\usepackage{amssymb}
\usepackage{amsmath}
\usepackage{graphicx}
\usepackage{latexsym}
\usepackage{amsfonts}

\numberwithin{equation}{section}

\newtheorem{theo}{Theorem}[section]

\newtheorem{pro}{Proposition}[section]
\newtheorem{cor}{Corollary}[section]
\newtheorem{exa}{Example}[section]
\newtheorem{axiom}{Definition}[section]
\newtheorem{rem}{Remark}[section]
\newtheorem{lem}{Lemma}[section]

\input{tcilatex}
\begin{document}
\title{ Stochastic Optimal Multi-Modes Switching with a \\Viscosity Solution Approach \thanks{The research
leading to these results has received funding from the European
Community's FP 7 Programme under contract agreement
PITN-GA-2008-213841, Marie Curie ITN "Controlled Systems".} }
\author{Brahim EL ASRI \thanks{%
Institut f\"{u}r Stochastik Friedrich-Schiller-Universit\"{a}t Jena
Ernst-Abbe-Platz 2 07743 Jena, Germany; e-mail:
brahim.el-asri@uni-jena.de }}
\date{}
\maketitle

\begin{abstract}
We consider the problem of optimal multi-modes switching in finite
horizon, when the state of the system, including the switching cost
functions are arbitrary ($g_{ij}(t,x)\geq 0$). We show existence of
the optimal strategy, and give when the optimal strategy is finite
via a verification theorem. Finally, when the state of the system is
a markov process, we show that the vector of value functions of the
optimal problem is the unique viscosity solution to the system of
$m$ variational partial differential inequalities with
inter-connected obstacles.
\end{abstract}

\noindent {$\mathbf{Keywords}$.} Real options, Backward stochastic
differential equations, Snell envelope, Stopping times, Switching,
Viscosity solution of PDEs, Variational inequalities

\noindent \textbf{AMS Classification subjects.} 60G40, 62P20, 91B99,
91B28, 35B37, 49L25

\medskip

\section{Introduction}

We consider a power plant which produces electricity and which has
several modes of production, e.g., the lower, the middle and the
intensive modes. The price $(X_{t})_{t\geq 0}$ of electricity in the
market fluctuates in reaction to many factors such as demand level,
weather conditions, unexpected outages etc. Moreover, electricity is
non-storable once produced, it should be almost immediately
consumed. Therefore, as a consequence, the station produces
electricity in its instantaneous most profitable mode known
that when the plant is in mode $i\in \mathcal{I}$, the yield per unit time $%
dt$ is given by means of $\psi _{i}(t,X_{t})dt$ and, on the other
hand, switching the plant from the mode $i$ to the mode $j$ is not
free and generates expenditures given by $g_{ij}(t,X_{t})$ and
possibly by other factors in the energy market.

The switching from one regime to another one is realized
sequentially at random times which are part of the decisions. So the
manager of the power plant faces two main issues:

$(i)$ when should she decide to switch the production from its
current mode to another one?

$(ii)$ to which mode the production has to be switched when the
decision of switching is made?

The manager faces the issue of finding the optimal strategy of
management of the plant. This is related with the price of the power
plant in the energy market.

Optimal switching problems for stochastic systems were studied by
several authors (see e.g. \cite{BO1, BS, CL, DP, DH, DHP, DZ2, E,
EH, HJ, TY, dz} and the references therein). The motivations are
mainly related to decision making in the economic sphere. In order
to tackle those problems, authors use mainly two approaches. Either
a probabilistic one \cite{DH, DHP, HJ} or an approach which uses
partial differential inequalities (PDIs for short)
\cite{BO1,CL,DZ2,EH, dz, TY}.

In the finite horizon framework Djehiche \textit{et al.} \cite{DHP}
have studied the multi-modes switching problem in using
probabilistic tools. They have proved existence of a solution and
found an optimal strategy when the
switching costs from state $i$ to state $j$ is strictly non--negative ($%
g_{ij}>\alpha >0$). The partial differential equation approach of
this work has been carried out by El Asri and Hamad\`{e}ne
\cite{EH}. We showed that when the price process $(X_{t}:t\geq 0)$
is solution of a Markovian stochastic differential equation, then
this problem is associated to a system of variational inequalities
with interconnected obstacles for which we provided a solution in
viscosity sense. This solution is bind to the value function of the
problem. Moreover the solution of the system is unique.

Using purely probabilistic tools such as the system of backward
stochastic
differential equations with oblique reflections (RBSDEs for short), Hamad%
\`{e}ne and Zhang \cite{HZ} have considered this optimal switching
problem
when the switching costs from state $i$ to state $j$ is non--negative $%
g_{ij} $. But in general case the optimal strategy may not exist.

The purpose of this work is to fill in this gap by providing a
solution to the optimal multiple switching problem using
probabilistic tools and partial differential equation approach.

We prove existence and provide a characterization of an optimal
strategy of
this problem when the payoff rates $\psi _{i}$ and the switching costs $%
g_{ij}\geq 0$ are adapted only to the filtration generated by a
Brownian motion. Later on, in the case when $X$ is a solution of a
SDE, we show that the value function of the problem is associated an
uplet of deterministic functions $(v^{1},\dots ,v^{m})$ which is the
unique solution of the following system of PDIs:
\begin{equation}
\left\{
\begin{array}{l}
\min \left\{ v_{i}(t,x)-\max\limits_{j\in \mathcal{I}^{-i}}\left\{
-g_{ij}(t,x)+v_{j}(t,x)\right\} ,-\partial _{t}v_{i}(t,x)-\mathcal{A}%
v_{i}(t,x)-\psi _{i}(t,x)\right\} =0\\
\forall \,\,(t,x)\in \lbrack 0,T]\times \mathbb{R}^{k},\,\,i\in \mathcal{I}%
=\{1,...,m\},\text{\thinspace\ and \thinspace\ }v_{i}(T,x)=0,%
\end{array}%
\right.  \label{sysintro}
\end{equation}%
where $\mathcal{A}$ an operator associated with a diffusion process and $%
\mathcal{I}^{-i}:=\mathcal{I}\setminus \{i\}$. It turns out that
this system is the deterministic version of the Verification Theorem
of the optimal multi-modes switching problem in infinite horizon.

This paper is organized as follows: In Section 2, we formulate the
problem and give the related definitions. In Section 3, we shall
introduce the optimal switching problem under consideration and give
its probabilistic Verification Theorem. It is expressed by means of
a Snell envelope. Then we introduce the approximating scheme which
enables us to construct a solution for the Verification Theorem.
Moreover we give some properties of that solution, especially the
dynamic programming principle. Section 4 is devoted to the
connection between the optimal switching problem, the Verification
Theorem and the associated system of PDIs. This connection is made
through BSDEs with one reflecting obstacle in the Markovian case.
Further we provide some estimate for the optimal strategy of the
switching problem which, in combination with the dynamic programming
principle, plays a crucial role in the proof of existence of a
solution for (\ref{sysintro}). In Section 5, we show that the
solution of PDIs is unique in the class of continuous functions
which satisfy a polynomial growth condition. In section 6 some
numerical examples are given. We close this paper an appendix in
which some technical results are proved.

\section{Assumptions and formulation of the problem}

Throughout this paper $T$ (resp. $k,\, d$) is a fixed real (resp.
integers) positive numbers. \newline Let \medskip

\noindent $(i)$ $b:[0,T]\times \mathbb{R}^{k}\rightarrow
\mathbb{R}^{k}$ and $\sigma :[0,T]\times \mathbb{R}^{k}\rightarrow
\mathbb{R}^{k\times d}$ be two continuous functions for which there
exists a constant $C>0$ such that for any $t\in \lbrack 0,T]$ and
$x,x^{\prime }\in \mathbb{R}^{k}$
\begin{equation}
|\sigma (t,x)-\sigma (t,x^{\prime })|+|b(t,x)-b(t,x^{\prime })|\leq
C|x-x^{\prime }|,  \label{regbs1}
\end{equation}%
$(ii)$ for $i,j\in \mathcal{I}=\{1,\ldots ,m\}$, $g_{ij}:[0,T]\times \mathbb{%
R}^{k}\rightarrow \mathbb{R}$ and $\psi _{i}:[0,T]\times \mathbb{R}%
^{k}\rightarrow \mathbb{R}$ are continuous functions and of
polynomial growth, $i.e.$ there exist some positive constants $C$
and $\gamma $ such that for each $i,j\in \mathcal{I}$:
\begin{equation}
|\psi _{i}(t,x)|+|g_{ij}(t,x)|\leq C(1+|x|^{\gamma }),\,\,\forall
\;(t,x)\in \lbrack 0,T]\times \mathbb{R}^{k},  \label{polycond}
\end{equation}%
$(iii)$ for any $(t,x)\in \lbrack 0,T]\times \mathbb{R}^{k}$,
$g_{ij}(t,x)$ are satisfying
\begin{equation}
g_{ii}(t,x)=0,\quad g_{ij}(t,x)\geq 0\quad \mbox{and}\quad
g_{ij}(t,x)+g_{jk}(t,x)>g_{ik}(t,x),\quad j\neq i,\ k\in
\mathcal{I},
\end{equation}%
which means that it is less expensive to switch directly in one step
from regime $i$ to $k$ than in two steps via an intermediate regime
$j$.

Moreover we assume that there exists a constant $\alpha >0$ such
that for any $(t,x)\in [0,T]\times \mathbb{R}^k$,
\begin{equation}
g_{ij}(t,x)+g_{ji}(t,x) > \alpha,\qquad i \neq j \in \mathcal{I} .
\end{equation}
This condition means that switching back and forth is not free.

We now consider the following system of $m$ variational inequalities
with inter-connected obstacles: $\forall \,\,i\in \mathcal{I}$
\begin{equation}
\left\{
\begin{array}{l}
\min \left\{ v_{i}(t,x)-\max\limits_{j\in \mathcal{I}^{-i}}\left\{
-g_{ij}(t,x)+v_{j}(t,x)\right\} ,-\partial _{t}v_{i}(t,x)-\mathcal{A}%
v_{i}(t,x)-\psi _{i}(t,x)\right\} =0, \\
v_{i}(T,x)=0,%
\end{array}%
\right.  \label{sysvi1}
\end{equation}%
where $\mathcal{A}$ is given by:
\begin{equation}
\mathcal{A}=\frac{1}{2}\sum_{i,j=1}^{m}(\sigma \sigma ^{\ast })_{ij}(t,x)%
\frac{\partial ^{2}}{\partial x_{i}\partial x_{j}}+\sum_{i=1}^{m}b_{i}(t,x)%
\frac{\partial }{\partial x_{i}}\,;  \label{generateur}
\end{equation}%
hereafter the superscript $(^{\ast })$ stands for the transpose,
$Tr$ is the trace operator and finally $\left\langle
x,y\right\rangle $ is the inner product of $x,y\in \mathbb{R}^{k}$.

The main objective of this paper is to focus on the uniqueness of
the solution in viscosity sense of (\ref{sysvi1}) whose definition
is:

\begin{axiom}
Let $(v_1,\ldots ,v_m)$ be a $m$-uplet of continuous functions defined on $%
[0,T]\times \mathbb{R}^k$, $\mathbb{R}$-valued and such that
$v_i(T,x)=0$
for any $x\in \mathbb{R}^k$ and $i\in \mathcal{I}$. The $m$-uplet $%
(v_1,\ldots ,v_m)$ is called:

\begin{itemize}
\item[$(i)$] a viscosity supersolution (resp. subsolution) of the system (%
\ref{sysvi1}) if for each fixed $i\in \mathcal{I}$, for any $%
(t_{0},x_{0})\in \lbrack 0,T]\times \mathbb{R}^{k}$ and any function $%
\varphi _{i}\in \mathcal{C}^{1,2}([0,T]\times \mathbb{R}^{k})$ such that $%
\varphi _{i}(t_{0},x_{0})=v_{i}(t_{0},x_{0})$ and $(t_{0},x_{0})$ is
a local maximum of $\varphi _{i}-v_{i}$ (resp. minimum), we have:
\begin{equation}
\min \left\{ v_{i}(t_{0},x_{0})-\max\limits_{j\in
\mathcal{I}^{-i}}\left\{
-g_{ij}(t_{0},x_{0})+v_{j}(t_{0},x_{0})\right\} ,-\partial
_{t}\varphi _{i}(t_{0},x_{0})-\mathcal{A}\varphi
_{i}(t_{0},x_{0})-\psi _{i}(t_{0},x_{0})\right\} \geq 0\,\,
(\mbox{resp.}\leq 0).
\end{equation}

\item[$(ii)$] a viscosity solution if it is both a viscosity supersolution
and subsolution. $\Box$
\end{itemize}
\end{axiom}

There is an equivalent formulation of this definition (see e.g.
\cite{CIL}) which we give because it will be useful later. So
firstly we define the notions of superjet and subjet of a continuous
function $v$.

\begin{axiom}
Let $v\in \mathcal{C}((0,T)\times \mathbb{R}^{k})$, $(t,x)$ an element of $%
(0,T)\times \mathbb{R}^{k}$ and finally $\mathbb{S}_{k}$ the set of
$k\times k$ symmetric matrices. We denote by $J^{2,+}v(t,x)$ (resp.
$J^{2,-}v(t,x)$),
the superjets (resp. the subjets) of $v$ at $(t,x)$, the set of triples $%
(p,q,X)\in \mathbb{R}\times \mathbb{R}^{k}\times \mathbb{S}_{k}$
such that:
\begin{eqnarray*}
v(s,y) &\leq &v(t,x)+p(s-t)+\langle q,y-x\rangle +\frac{1}{2}\langle
X(y-x),y-x\rangle +o(|s-t|+|y-x|^{2}) \\
&& \\
(resp.\quad v(s,y) &\geq &v(t,x)+p(s-t)+\langle q,y-x\rangle +\frac{1}{2}%
\langle X(y-x),y-x\rangle +o(|s-t|+|y-x|^{2})).
\end{eqnarray*}
\end{axiom}

Note that if $\varphi -v$ has a local maximum (resp. minimum) at
$(t,x)$, then we obviously have:
\begin{equation*}
\left( D_{t}\varphi (t,x),D_{x}\varphi (t,x),D_{xx}^{2}\varphi
(t,x)\right) \in J^{2,-}v(t,x)\,\,\,(\mbox{resp.
}J^{2,+}v(t,x)).\Box
\end{equation*}

We now give an equivalent definition of a viscosity solution of the
parabolic system with inter-connected obstacles (\ref{sysvi1}).

\begin{axiom}
Let $(v_1,\ldots ,v_m)$ be a $m$-uplet of continuous functions defined on $%
[0,T]\times \mathbb{R}^k$, $\mathbb{R}$-valued and such that
$(v_1,\ldots ,v_m)(T,x)=0$ for any $x\in \mathbb{R}^k$. The
$m$-uplet $(v_1,\ldots ,v_m)$ is called a viscosity supersolution
(resp. subsolution) of (\ref{sysvi1}) if
for any $i\in \mathcal{I}$, $(t,x)\in (0,T)\times \mathbb{R}^k$ and $%
(p,q,X)\in J^{2,-} v_i (t,x)$ (resp. $J^{2,+} v_i (t,x)$),
\begin{equation*}
min \left\{v_i(t,x)-\max\limits_{j\in\mathcal{\
I}^{-i}}(-g_{ij}(t,x)+ v_j(t,x)),-p -\frac{1}{2}Tr[\sigma^{*} X
\sigma] -\langle b,q \rangle-\psi_{i}(t,x)\right\}\geq 0 \,\,(resp.
\leq 0).
\end{equation*}
It is called a viscosity solution it is both a viscosity subsolution
and supersolution. $\Box$
\end{axiom}

As pointed out previously we will show that system (\ref{sysvi1})
has a unique solution in viscosity sense. This system is the
deterministic version of the optimal $m$-states switching problem
will describe briefly in the next section.

\section{The optimal $m$-states switching problem}

\subsection{Setting of the problem}

Let $(\Omega, \mathcal{F}, P)$ be a fixed probability space on which
is defined a standard $d$-dimensional Brownian motion
$B=(B_t)_{0\leq t\leq T}$ whose natural filtration is
$(\mathcal{F}_t^0:=\sigma \{B_s, s\leq t\})_{0\leq t\leq T}$. Let
$\mathbf{F}=(\mathcal{F}_t)_{0\leq t\leq T}$ be
the completed filtration of $(\mathcal{F}_t^0)_{0\leq t\leq T}$ with the $P$%
-null sets of $\mathcal{F}$. \newline Let:

- $\mathcal{P}$ be the $\sigma$-algebra on $[0,T]\times \Omega$ of $\mathbf{F%
}$-progressively measurable sets;

- $\mathcal{M}^{2,k}$ be the set of $\mathcal{P}$-measurable and $\mathbb{R}%
^k$-valued processes $w=(w_t)_{t\leq T}$ such that $E[\int_0^T|w_s|^2ds]<%
\infty$ and $\mathcal{S}^2$ be the set of $\mathcal{P}$-measurable,
continuous processes ${w}=({w}_t)_{t\leq T}$ such that $E[\sup_{t\leq T}|{w}%
_t|^2]<\infty$;

- for any stopping time $\tau \in [0,T]$, $\mathcal{T}_\tau$ denotes
the set of all stopping times $\theta$ such that $\tau \leq \theta
\leq T$. \medskip

Let $\mathcal{I}$ be the set of all possible activity modes of the
production of a power plant. A management strategy of the plant
consists, on
the one hand, of the choice of a sequence of nondecreasing stopping times $%
(\tau _{n})_{n\geq 1}$ (i.e. $\tau _{n}\leq \tau _{n+1}$ and $\tau
_{0}=0$) where the manager decides to switch the activity from its
current mode to another one. On the other hand, it consists of the
choice of the mode $\xi _{n}$, which is an $\mathcal{F}_{\tau
_{n}}$-measurable random variable taking values in $\mathcal{I}$, to
which the production is switched at $\tau _{n}$ from its current
mode. Therefore the admissible management strategies of the plant
are the pairs $(\delta ,\xi ):=((\tau _{n})_{n\geq 1},(\xi
_{n})_{n\geq 1})$ and the set of these strategies is denoted by $\mathcal{D}$%
.

Let $X:=(X_{t})_{0\leq t\leq T}$ be an $\mathcal{P}$-measurable, $\mathbb{R}%
^{k}$-valued continuous stochastic process which stands for the
market price of $k$ factors which determine the market price of the
commodity. Assuming
that the production activity is in mode 1 at the initial time $t=0$, let $%
(u_{t})_{t\leq T}$ denote the indicator of the production activity's
mode at time $t\in \lbrack 0,T]$ :
\begin{equation}
u_{t}=\mathbf{1}_{[0,\tau _{1}]}(t)+\sum_{n\geq 1}\xi
_{n}\mathbf{1}_{(\tau _{n},\tau _{n+1}]}(t).
\end{equation}%
Then for any $t\leq T$, the state of the whole economic system
related to the project at time $t$ is represented by the vector:
\begin{equation}
(t,X_{t},u_{t})\in \lbrack 0,T]\times \mathbb{R}^{k}\times
\mathcal{I}.
\end{equation}%
Finally, let $\psi _{i}(t,X_{t})$ be the instantaneous profit when
the system is in state $(t,X_{t},i)$, and for $i,j\in \mathcal{I}$
$\ i\neq j$, let $g_{ij}(t,X_{t})$ denote the switching cost of the
production at time $t$ from current mode $i$ to another mode $j$.
Then if the plant is run under the strategy $(\delta ,\xi )=((\tau
_{n})_{n\geq 1},(\xi _{n})_{n\geq 1})$ the expected total profit is
given by:
\begin{equation*}
J(\delta ,\xi )=E\left[ \int_{0}^{T}\psi
_{u_{s}}(s,X_{s})ds-\sum_{n\geq 1}g_{u_{\tau _{n-1}}u_{\tau
_{n}}}(\tau _{n},X_{\tau _{n}})\mathbf{1}_{[\tau _{n}<T]}\right] .
\end{equation*}%
Therefore the problem we are interested in is to find an optimal strategy $%
i.e.$ a strategy $(\delta ^{\ast },\xi ^{\ast })$ such that
$J(\delta ^{\ast },\xi ^{\ast })\geq J(\delta ,\xi )$ for any
$(\delta ,\xi )\in \mathcal{D}$.

Note that in order that the quantity $J(\delta ,\xi )$ makes sense,
we assume throughout this paper that, for any $i,j\in \mathcal{I}$
the processes $(\psi _{i}(t,X_{t}))_{t\leq T}$ and
$(g_{ij}(t,X_{t}))_{t\leq T}$ belong to $\mathcal{M}^{2,1}$ and
$\mathcal{S}^{2}$ respectively. There is
one to one correspondence between the pairs $(\delta ,\xi )$ and the pairs $%
(\delta ,u)$. Therefore throughout this paper one refers indifferently to $%
(\delta ,\xi )$ or $(\delta ,u)$.

\subsection{The Verification Theorem}

To tackle the problem described above Djehiche \textit{et al.}
\cite{DHP} have introduced a Verification Theorem which is expressed
by means of Snell envelope of processes. The Snell envelope of a
stochastic process $(\eta _{t})_{t\leq T}$ of $\mathcal{S}^{2}$
(with a possible positive jump at $T$)
is the lowest supermartingale $R(\eta ):=(R(\eta )_{t})_{t\leq T}$ of $%
\mathcal{S}^{2}$ such that for any $t\leq T$, $R(\eta )_{t}\geq \eta
_{t}$. It has the following expression:
\begin{equation*}
\forall \;t\leq T,\ R(\eta )_{t}=\limfunc{esssup}_{\tau \in \mathcal{T}%
_{t}}E[\eta _{\tau }|\mathcal{F}_{t}]\mbox{ and satisfies }R(\eta
)_{T}=\eta _{T}.
\end{equation*}%
For more details owe refer to \cite{CK, Elka, ham}. \medskip

The Verification Theorem for the $m$-states optimal switching
problem is the following:

\begin{theo}
\label{thmverif} Assume that there exist $m$ processes $%
(Y^{i}:=(Y_{t}^{i})_{0\leq t\leq T},i=1,\ldots ,m)$ of
$\mathcal{S}^{2}$ such that:
\begin{equation}
\forall \;t\leq T,\,\,Y_{t}^{i}=\limfunc{esssup}_{\tau \geq
t}E\left[
\int_{t}^{{\tau }}\psi _{i}(s,X_{s})ds+\max\limits_{j\in \mathcal{I}%
^{-i}}(-g_{ij}(\tau ,X_{\tau })+Y_{\tau }^{j})\mathbf{1}_{[\tau <T]}|%
\mathcal{F}_{t}\right] ,\;Y_{T}^{i}=0.  \label{eqvt}
\end{equation}%
Then:

\begin{itemize}
\item[(i)] $Y_{0}^{1}=\sup\limits_{(\delta ,\xi )\in \mathcal{D}}J(\delta
,u).$

\item[(ii)] Define the sequence of $\mathbf{F}$-stopping times $\delta
^{\ast }=(\tau _{n}^{\ast })_{n\geq 1}$ as follows :
\begin{equation*}
\begin{array}{lll}
\tau _{1}^{\ast } & = & \inf \left\{ s\geq 0,\quad
Y_{s}^{1}=\max\limits_{j\in {\mathcal{\ I}^{-1}}%
}(-g_{1j}(s,X_{s})+Y_{s}^{j})\right\} \wedge T, \\
\tau _{n}^{\ast } & = & \inf \left\{ s\geq \tau _{n-1}^{\ast },\quad
Y_{s}^{u_{\tau _{n-1}^{\ast }}}=\max\limits_{k\in
\mathcal{I}\backslash \{u_{\tau _{n-1}^{\ast }}\}}(-g_{u_{\tau
_{n-1}^{\ast
}}k}(s,X_{s})+Y_{s}^{k})\right\} \wedge T,\quad \mbox{for}\quad n\geq 2,%
\end{array}%
\end{equation*}%
where:

\begin{itemize}
\item[$\bullet $] $u_{\tau _{1}^{\ast }}=\sum\limits_{j\in \mathcal{\ I}}j%
\mathbf{1}_{\{\max\limits_{k\in \mathcal{\ I}^{-1}}(-g_{1k}(\tau
_{1}^{\ast },X_{\tau _{1}^{\ast }})+Y_{\tau _{1}^{\ast
}}^{k})=-g_{1j}(\tau _{1}^{\ast },X_{\tau _{1}^{\ast }})+Y_{\tau
_{1}^{\ast }}^{j}\}};$

\item[$\bullet $] for any $n\geq 1$ and $t\geq \tau _{n}^{\ast },$ $%
Y_{t}^{u_{\tau _{n}^{\ast }}}=\sum\limits_{j\in \mathcal{I}}\mathbf{1}%
_{[u_{\tau _{n}^{\ast }}=j]}Y_{t}^{j}$

\item[$\bullet $] for any $n\geq 2,\,\,u_{\tau _{n}^{\ast }}=l$ on the set
\begin{equation*}
\left\{ \max\limits_{k\in \mathcal{I}\backslash \{{u_{\tau _{n-1}^{\ast }}}%
\}}(-g_{u_{\tau _{n-1}^{\ast }}k}(\tau _{n}^{\ast },X_{\tau
_{n}^{\ast }})+Y_{\tau _{n}^{\ast }}^{k})=-g_{u_{\tau _{n-1}^{\ast
}l}}(\tau _{n}^{\ast },X_{\tau _{n}^{\ast }})+Y_{\tau _{n}^{\ast
}}^{l}\right\}
\end{equation*}%
with $g_{u_{\tau _{n-1}^{\ast }k}}(\tau _{n}^{\ast },X_{\tau
_{n}^{\ast }})=\sum\limits_{j\in \mathcal{I}}\mathbf{1}_{[u_{\tau
_{n-1}^{\ast
}}=j]}g_{jk}(\tau _{n}^{\ast },X_{\tau _{n}^{\ast }})$ and $\mathcal{I}%
\backslash \{u_{\tau _{n-1}^{\ast }}\}=\sum\limits_{j\in \mathcal{I}}\mathbf{%
1}_{[u_{\tau _{n-1}^{\ast }}=j]}\mathcal{I}^{-j}$.\newline Then the
strategy $(\delta ^{\ast },u^{\ast })$ satisfies
\begin{equation*}
E\left[ \sum_{k\geq 1}g_{u_{\tau _{k-1}^{\ast }}u_{\tau _{k}^{\ast
}}}({\tau
_{k}^{\ast }},X_{{\tau _{k}^{\ast }}})\mathbf{1}_{[\tau _{k}^{\ast }<T]}%
\right] <+\infty
\end{equation*}%
and it is optimal.

\item[(iii)] If
\begin{equation*}
E\left[ \sum\limits_{k\geq 1}g_{u_{\tau _{k}^{\ast }}u_{\tau
_{k-1}^{\ast }}}({\tau _{k}^{\ast }},X_{{\tau _{k}^{\ast
}}})\mathbf{1}_{[\tau _{k}^{\ast }<T]}\right] <+\infty
\end{equation*}%
or $g_{ij}$ is constant, then the optimal strategy $(\delta ^{\ast
},u^{\ast })$ is finite.
\end{itemize}
\end{itemize}
\end{theo}
$Proof .$ The proof is divided in four steps\\ \textbf{Step 1.
}(i)\textbf{\ }It consists in showing that for any $t\leq T,$
$Y_{t}^{i}$, as defined by (\ref{eqvt}), is the expected total
profit or the value function of the optimal problem, given that the
system is in mode $i$ at time $t$. More precisely,
\begin{equation*}
Y_{t}^{i}=\limfunc{esssup}_{(\delta ,u)\in \mathcal{D}_{t}}E\left[
\int_{t}^{T}\psi _{i}(s,X_{s})ds-\sum_{k\geq 1}g_{u_{\tau
_{k-1}}u_{\tau
_{k}}}(\tau _{k},X_{{\tau _{k}}})\mathbf{1}_{[\tau _{k}<T]}|\mathcal{F}_{t}%
\right] ,
\end{equation*}

where $\mathcal{D}_{t}$ is the set of strategies such that $\tau _{1}\geq t$%
, $P$-a.s. if at time $t$ the system is in the mode $i$.

Let us admit for a moment the following Lemma whose proof is given
in the appendix.

\begin{lem}
\label{env-sne} For every $\tau _{1}^{\ast }\leq t\leq T$.
\begin{equation}
Y_{t}^{u_{\tau _{1}^{\ast }}}=\limfunc{esssup}_{\tau \geq t}E\left[
\int_{t}^{\tau }\psi _{u_{\tau _{1}^{\ast
}}}(s,X_{s})ds+\max\limits_{j\in \mathcal{I}^{-u_{\tau _{1}^{\ast
}}}}(-g_{u_{\tau _{1}^{\ast }}j}(\tau ,X_{\tau })+Y_{\tau
}^{j})\mathbf{1}_{[\tau <T]}|\mathcal{F}_{t}\right] .\;\Box
\label{env1}
\end{equation}
\end{lem}

 From properties of the Snell envelope and at time $t=0$
the system is in mode $1$, we have:
\begin{eqnarray*}
Y_{0}^{1} &=&E\left[ \int_{0}^{\tau _{1}^{\ast }}\psi
_{1}(s,X_{s})ds+\max\limits_{j\in \mathcal{I}^{-1}}(-g_{1j}({\tau
_{1}^{\ast }},X_{\tau _{1}^{\ast }})+Y_{\tau _{1}^{\ast
}}^{j})\mathbf{1}_{[\tau
_{1}^{\ast }<T]}\right]  \\
&& \\
&=&E\left[ \int_{0}^{\tau _{1}^{\ast }}\psi
_{1}(s,X_{s})ds+(-g_{1u_{\tau _{1}^{\ast }}}({\tau _{1}^{\ast
}},X_{\tau _{1}^{\ast }})+Y_{\tau _{1}^{\ast }}^{u_{\tau _{1}^{\ast
}}})\mathbf{1}_{[\tau _{1}^{\ast }<T]}\right] .
\end{eqnarray*}%
Now, from Lemma \ref{env-sne} and the definition of $\tau _{2}^{\ast
}$ we have:
\begin{eqnarray*}
Y_{\tau _{1}^{\ast }}^{u_{\tau _{1}^{\ast }}} &=&E\left[ \int_{\tau
_{1}^{\ast }}^{\tau _{2}^{\ast }}\psi _{u_{\tau _{1}^{\ast
}}}(s,X_{s})ds+\max\limits_{j\in \mathcal{I}^{-u_{\tau _{1}^{\ast
}}}}(-g_{u_{\tau _{1}^{\ast }}j}(\tau _{2}^{\ast },X_{\tau
_{2}^{\ast
}})+Y_{\tau _{2}^{\ast }}^{j})\mathbf{1}_{[\tau _{2}^{\ast }<T]}|\mathcal{F}%
_{\tau _{1}^{\ast }}\right]  \\
&& \\
&=&E\left[ \int_{\tau _{1}^{\ast }}^{\tau _{2}^{\ast }}\psi
_{u_{\tau _{1}^{\ast }}}(s,X_{s})ds+(-g_{u_{\tau _{1}^{\ast
}}u_{\tau _{2}^{\ast }}}(\tau _{2}^{\ast },X_{\tau _{2}^{\ast
}})+Y_{\tau _{2}^{\ast }}^{u_{\tau _{2}^{\ast }}})\mathbf{1}_{[\tau
_{2}^{\ast }<T]}|\mathcal{F}_{\tau _{1}^{\ast }}\right] .
\end{eqnarray*}

It implies that
\begin{eqnarray*}
Y_{0}^{1} &=&E\left[ \int_{0}^{\tau _{1}^{\ast }}\psi
_{1}(s,X_{s})ds-g_{1u_{\tau _{1}^{\ast }}}({\tau _{1}^{\ast
}},X_{\tau
_{1}^{\ast }})\mathbf{1}_{[\tau _{1}^{\ast }<T]}\right]  \\
&& \\
&&+E\left[ \int_{\tau _{1}^{\ast }}^{\tau _{2}^{\ast }}\psi
_{u_{\tau
_{1}^{\ast }}}(s,X_{s})ds+(-g_{u_{\tau _{1}^{\ast }}u_{\tau _{2}^{\ast }}}({%
\tau _{2}^{\ast }},X_{\tau _{2}^{\ast }})+Y_{\tau _{2}^{\ast
}}^{u_{\tau _{2}^{\ast }}})\mathbf{1}_{[\tau _{2}^{\ast
}<T]}|\mathcal{F}_{\tau
_{1}^{\ast }}\right]  \\
&& \\
&=&E\left[ \int_{0}^{\tau _{1}^{\ast }}\psi
_{1}(s,X_{s})ds+\int_{\tau _{1}^{\ast }}^{\tau _{2}^{\ast }}\psi
_{u_{\tau _{1}^{\ast
}}}(s,X_{s})ds\right.  \\
&& \\
&&\left. -g_{1u_{\tau _{1}^{\ast }}}({\tau _{1}^{\ast }},X_{\tau
_{1}^{\ast }})\mathbf{1}_{[\tau _{1}^{\ast }<T]}-g_{u_{\tau
_{1}^{\ast }}u_{\tau _{2}^{\ast }}}({\tau _{2}^{\ast }},X_{\tau
_{2}^{\ast }})\mathbf{1}_{[\tau
_{2}^{\ast }<T]}+Y_{\tau _{2}^{\ast }}^{u_{\tau _{2}^{\ast }}}\mathbf{1}%
_{[\tau _{2}^{\ast }<T]}\right] ,
\end{eqnarray*}%
since $[\tau _{2}^{\ast }<T]\subset \lbrack \tau _{1}^{\ast }<T]$.
Therefore
\begin{equation}
Y_{0}^{1}=E\left[ \int_{0}^{\tau _{2}^{\ast }}\psi
(s,X_{s},u_{s})ds-g_{1u_{\tau _{1}^{\ast }}}({\tau _{1}^{\ast
}},X_{\tau _{1}^{\ast }})\mathbf{1}_{[\tau _{1}^{\ast
}<T]}-g_{u_{\tau _{1}^{\ast
}}u_{\tau _{2}^{\ast }}}({\tau _{2}^{\ast }},X_{\tau _{2}^{\ast }})\mathbf{1}%
_{[\tau _{2}^{\ast }<T]}+Y_{\tau _{2}^{\ast }}^{u_{\tau _{2}^{\ast }}}%
\mathbf{1}_{[\tau _{2}^{\ast }<T]}\right] ,  \label{principe}
\end{equation}%
since between $0$ and $\tau _{1}^{\ast }$ (resp. $\tau _{1}^{\ast }$ and $%
\tau _{2}^{\ast }$) the production is in regime $1$ (resp. regime
$u_{\tau _{1}^{\ast }}$) and then $u_{t}=1$ (resp. $u_{t}=u_{\tau
_{1}^{\ast }}$) which implies that
\begin{equation*}
\int_{0}^{\tau _{2}^{\ast }}\psi (s,X_{s},u_{s})ds=\int_{0}^{\tau
_{1}^{\ast }}\psi _{1}(s,X_{s})ds+\int_{\tau _{1}^{\ast }}^{\tau
_{2}^{\ast }}\psi _{u_{\tau _{1}^{\ast }}}(s,X_{s})ds.
\end{equation*}%
Repeating this reasoning as many times as necessary we obtain that for any $%
n\geq 0,$
\begin{equation*}
Y_{0}^{1}=E\left[ \int_{0}^{\tau _{n}^{\ast }}\psi
(s,X_{s},u_{s})ds-\sum_{k=1}^{n}g_{u_{\tau _{k-1}^{\ast }}u_{\tau
_{k}^{\ast }}}({\tau _{k}^{\ast }},X_{{\tau _{k}^{\ast
}}})\mathbf{1}_{[\tau _{k}^{\ast }<T]}+Y_{\tau _{n}^{\ast
}}^{u_{\tau _{n}^{\ast }}}\mathbf{1}_{[\tau _{n}^{\ast }<T]}\right]
.
\end{equation*}%
Then, the strategy $(\delta ^{\ast },u^{\ast })$ satisfies
\begin{equation*}
E\left[ \sum\limits_{k\geq 1}g_{u_{\tau _{k-1}^{\ast }}u_{\tau
_{k}^{\ast }}}({\tau _{k}^{\ast }},X_{{\tau _{k}^{\ast
}}})\mathbf{1}_{[\tau _{k}^{\ast }<T]}\right] <+\infty .
\end{equation*}%
If not $Y_{0}^{1}=-\infty $ which contradicts the assumption
$Y^{i}\in \mathcal{S}^{2}$. Therefore, taking the limit as
$n\rightarrow +\infty $ we obtain $Y_{0}^{1}=J(\delta ^{\ast
},u^{\ast })$. \\ \textbf{Step 2.} (ii) We show that the strategy
$(\delta ^{\ast },u^{\ast })$
it is optimal i.e. $J(\delta ^{\ast },u^{\ast })\geq J(\delta ,u)$ for any $%
(\delta ,u)\in \mathcal{D}$.\newline The definition of the Snell
envelope yields
\begin{eqnarray*}
Y_{0}^{1} &\geq &E\left[ \int_{0}^{\tau _{1}}\psi
_{1}(s,X_{s})ds+\max\limits_{j\in \mathcal{I}^{-1}}(-g_{1j}(\tau
_{1},X_{\tau _{1}})+Y_{\tau _{1}}^{j})\mathbf{1}_{[\tau _{1}<T]}\right]  \\
&& \\
&\geq &E\left[ \int_{0}^{\tau _{1}}\psi
_{1}(s,X_{s})ds+(-g_{1u_{\tau
_{1}^{\ast }}}(\tau _{1},X_{\tau _{1}})+Y_{\tau _{1}}^{u_{\tau _{1}}})%
\mathbf{1}_{[\tau _{1}<T]}\right] .
\end{eqnarray*}%
But, once more using a similar characterization as (\ref{env1}), we
get
\begin{eqnarray*}
Y_{\tau _{1}}^{u_{\tau _{1}}} &\geq &E\left[ \int_{\tau _{1}}^{\tau
_{2}}\psi _{u_{\tau _{1}}}(s,X_{s})ds+\max\limits_{j\in \mathcal{I}%
^{-u_{\tau _{1}}}}(-g_{u_{\tau _{1}}j}(\tau _{2},X_{\tau
_{2}})+Y_{\tau
_{2}}^{j})\mathbf{1}_{[\tau _{2}<T]}|\mathcal{F}_{\tau _{1}}\right]  \\
&& \\
&\geq &E\left[ \int_{\tau _{1}}^{\tau _{2}}\psi _{u_{\tau
_{1}}}(s,X_{s})ds+(-g_{u_{\tau _{1}}u_{\tau _{2}}}(\tau _{2},X_{\tau
_{2}})+Y_{\tau _{2}}^{u_{\tau _{2}}})\mathbf{1}_{[\tau _{2}<T]}|\mathcal{F}%
_{\tau _{1}}\right] .
\end{eqnarray*}%
Therefore,
\begin{eqnarray*}
Y_{0}^{1} &\geq &E\left[ \int_{0}^{\tau _{1}}\psi
_{1}(s,X_{s})ds-g_{1u_{\tau _{1}}}(\tau _{1},X_{\tau _{1}})\right]  \\
&& \\
&&+E\left[ \int_{\tau _{1}}^{\tau _{2}}\psi _{u_{\tau
_{1}}}(s,X_{s})ds+(-g_{u_{\tau _{1}}u_{\tau _{2}}}(\tau _{2},X_{\tau
_{2}})+Y_{\tau _{2}}^{u_{\tau _{2}}})\mathbf{1}_{[\tau _{2}<T]}\right]  \\
&& \\
&=&E\left[ \int_{0}^{\tau _{2}}\psi _{u_{s}}(s,X_{s})ds-g_{1u_{\tau
_{1}}}(\tau _{1},X_{\tau _{1}})\mathbf{1}_{[\tau _{1}<T]}-g_{u_{\tau
_{1}}u_{\tau _{2}}}(\tau _{2},X_{\tau _{2}})+Y_{\tau _{2}}^{u_{\tau _{2}}}%
\mathbf{1}_{[\tau _{2}<T]}\right] .
\end{eqnarray*}%
Repeat this argument $n$ times to obtain
\begin{equation*}
Y_{0}^{1}\geq E\left[ \int_{0}^{\tau _{n}}\psi
_{u_{s}}(s,X_{s})ds-\sum\limits_{k=1}^{n}g_{u_{\tau _{k-1}}u_{\tau
_{k}}}(\tau _{k},X_{{\tau _{k}}})\mathbf{1}_{[\tau _{k}<T]}+Y_{\tau
_{n}}^{u_{\tau _{n}}}\mathbf{1}_{[\tau _{n}<T]}\right] .
\end{equation*}%
Finally, taking the limit as $n\rightarrow +\infty $ yields
\begin{equation*}
Y_{0}^{1}\geq E\left[ \int_{0}^{T}\psi
_{u_{s}}(s,X_{s})ds-\sum_{k\geq
1}g_{u_{\tau _{k-1}}u_{\tau _{k}}}(\tau _{k},X_{{\tau _{k}}})\mathbf{1}%
_{[\tau _{k}<T]}\right] .
\end{equation*}%
Hence, the strategy $(\delta ^{\ast },u^{\ast })$ is optimal.\\
\textbf{Step. 3} (iii) Next, we show that the strategy $(\tau
_{n}^{\ast })_{n\geq 1}$ is finite if
\begin{equation*}
E\left[ \sum\limits_{k\geq 1}g_{u_{\tau _{k}^{\ast }}u_{\tau
_{k-1}^{\ast }}}({\tau _{k}^{\ast }},X_{{\tau _{k}^{\ast
}}})\mathbf{1}_{[\tau _{k}^{\ast }<T]}\right] <+\infty .
\end{equation*}%
Indeed, let $A=\{\omega :\tau _{n}^{\ast }(\omega )<T,\,\forall \ n\geq 1\}$%
. If $P(A)>0$, then from (\ref{principe}) we have for any $n\geq 1$,
\begin{eqnarray*}
&&Y_{0}^{1}-E\left[ \sum\limits_{k\geq 1}g_{u_{\tau _{k}^{\ast
}}u_{\tau
_{k-1}^{\ast }}}({\tau _{k}^{\ast }},X_{{\tau _{k}^{\ast }}})\mathbf{1}%
_{[\tau _{k}^{\ast }<T]}\right]  \\
&& \\
&\leq &E\left[ \int_{0}^{\tau _{n}^{\ast }}\max\limits_{1\leq j\leq
m}\mid \psi _{j}(s,X_{s})\mid ds-\sum_{k=1}^{n}g_{u_{\tau
_{k-1}^{\ast }}u_{\tau
_{k}^{\ast }}}({\tau _{k}^{\ast }},X_{{\ \tau _{k}^{\ast }}})\mathbf{1}%
_{[\tau _{k}^{\ast }<T]}+Y_{\tau _{n}^{\ast }}^{u_{\tau _{n}^{\ast }}}%
\mathbf{1}_{[\tau _{n}^{\ast }<T]}\right]  \\
&& \\
&&-E\left[ \sum\limits_{k=1}^{n}g_{u_{\tau _{k}^{\ast }}u_{\tau
_{k-1}^{\ast }}}({\tau _{k}^{\ast }},X_{{\tau _{k}^{\ast
}}})\mathbf{1}_{[\tau _{k}^{\ast
}<T]}\right]  \\
&& \\
&=&E\left[ \int_{0}^{\tau _{n}^{\ast }}\max\limits_{1\leq j\leq
m}\mid \psi _{j}(s,X_{s})\mid ds-\left( \sum\limits_{k=1}^{n}\left(
g_{u_{\tau _{k-1}^{\ast }}u_{\tau _{k}^{\ast }}}({\tau _{k}^{\ast
}},X_{{\tau _{k}^{\ast }}})+g_{u_{\tau _{k}^{\ast }}u_{\tau
_{k-1}^{\ast }}}({\tau _{k}^{\ast }},X_{{\tau _{k}^{\ast }}})\right)
\mathbf{1}_{[\tau _{k}^{\ast
}<T]}\right) \mathbf{1}_{A}\right.  \\
&& \\
&&\left. -\left( \sum\limits_{k=1}^{n}\left( g_{u_{\tau _{k-1}^{\ast
}}u_{\tau _{k}^{\ast }}}({\tau _{k}^{\ast }},X_{{\tau _{k}^{\ast }}%
})+g_{u_{\tau _{k}^{\ast }}u_{\tau _{k-1}^{\ast }}}({\tau _{k}^{\ast }},X_{{%
\tau _{k}^{\ast }}})\right) \mathbf{1}_{[\tau _{k}^{\ast
}<T]}\right)
\mathbf{1}_{\overline{A}}+Y_{\tau _{n}^{\ast }}^{u_{\tau _{n}^{\ast }}}%
\mathbf{1}_{[\tau _{n}^{\ast }<T]}\right]  \\
&& \\
&<&E\left[ \int_{0}^{\tau _{n}^{\ast }}\max\limits_{1\leq j\leq
m}\mid \psi _{j}(s,X_{s})\mid ds-\left( \sum\limits_{k=1}^{n}\alpha
\mathbf{1}_{[\tau _{k}^{\ast }<T]}\right) \mathbf{1}_{A}\,-\left(
\sum\limits_{k=1}^{n}\alpha \mathbf{1}_{[\tau _{k}^{\ast
}<T]}\right) \mathbf{1}_{\overline{A}}+Y_{\tau
_{n}^{\ast }}^{u_{\tau _{n}^{\ast }}}\mathbf{1}_{[\tau _{n}^{\ast }<T]}%
\right] ,
\end{eqnarray*}%
since $g_{ij}(t,x)+g_{ji}(t,x)>\alpha $. Then the right-hand side
converge
to $-\infty $ as $n\rightarrow \infty $. But this is contradictory because $%
Y^{i}$ belong to $\mathcal{S}^{2}$, $\psi _{i}(.,X)\in
\mathcal{M}^{2,1}$ and $E[\sum_{k\geq 1}g_{u_{\tau _{k}^{\ast
}}u_{\tau _{k-1}^{\ast }}}({\tau _{k}^{\ast }},X_{{\tau _{k}^{\ast
}}})\mathbf{1}_{[\tau _{k}^{\ast }<T]}]<+\infty $. Henceforth the
strategy is finite.
\\ \textbf{Step 4.} (iii) To complete the proof it remains to show that
the strategy $(\tau _{n}^{\ast })_{n\geq 1}$ is finite when $g_{ij}$
is constant. Indeed let $A=\{\omega ,\,\tau _{n}^{\ast }(\omega
)<T,\,\forall \;n\geq 1\}$. If $P(A)>0$, then from (\ref{principe})
we have for any $n\geq 1$,
\begin{equation*}
Y_{0}^{1}\leq E\left[ \int_{0}^{\tau _{nm}^{\ast
}}\max\limits_{1\leq j\leq m}\mid \psi _{j}(s,X_{s})\mid
ds-\sum\limits_{k=1}^{nm}g_{u_{\tau _{k-1}^{\ast }}u_{\tau
_{k}^{\ast }}}\mathbf{1}_{[\tau _{k}^{\ast }<T]}+Y_{\tau _{nm}^{\ast
}}^{u_{\tau _{nm}^{\ast }}}\mathbf{1}_{[\tau _{nm}^{\ast
}<T]}\right] .
\end{equation*}%
We show by induction on $n$ that for all $n\geq 1$,
\begin{equation}
-\sum_{1\leq k\leq nm}g_{u_{\tau _{k-1}^{\ast }}u_{\tau _{k}^{\ast }}}%
\mathbf{1}_{[\tau _{k}^{\ast }<T]}\leq -\alpha n\mathbf{1}_{[\tau
_{nm}^{\ast }<T]}.  \label{somme}
\end{equation}%
Indeed, the above assertion is obviously true for $n=1$. Suppose now
it holds true at step $n$. Then, at step $n+1$, we have
\begin{eqnarray*}
-\sum\limits_{k=1}^{(n+1)m}g_{u_{\tau _{k-1}^{\ast }}u_{\tau _{k}^{\ast }}}%
\mathbf{1}_{[\tau _{k}^{\ast }<T]}
&=&-\sum\limits_{k=1}^{nm}g_{u_{\tau _{k-1}^{\ast }}u_{\tau
_{k}^{\ast }}}\mathbf{1}_{[\tau _{k}^{\ast
}<T]}-\sum\limits_{k=nm+1}^{(n+1)m}g_{u_{\tau _{k-1}^{\ast }}u_{\tau
_{k}^{\ast }}}\mathbf{1}_{[\tau _{k}^{\ast }<T]} \\
&& \\
&\leq &-\alpha n\mathbf{1}_{[\tau _{nm}^{\ast }<T]}-\alpha
\mathbf{1}_{[\tau _{(n+1)m}^{\ast }<T]}\leq -\alpha
(n+1)\mathbf{1}_{[\tau _{(n+1)m}^{\ast }<T]}.
\end{eqnarray*}%
It follow that
\begin{eqnarray*}
Y_{0}^{1} &\leq &E\left[ \int_{0}^{\tau _{nm}^{\ast
}}\max\limits_{1\leq j\leq m}\mid \psi _{j}(s,X_{s})\mid ds-\alpha
n\mathbf{1}_{[\tau _{nm}^{\ast }<T]}+Y_{\tau _{nm}^{\ast }}^{u_{\tau
_{nm}^{\ast }}}\mathbf{1}_{[\tau
_{nm}^{\ast }<T]}\right]  \\
&& \\
&=&E\left[ \int_{0}^{\tau _{nm}^{\ast }}\max\limits_{1\leq j\leq
m}\mid \psi
_{j}(s,X_{s})\mid ds-\alpha n\mathbf{1}_{[\tau _{nm}^{\ast }<T]}\mathbf{1}%
_{A}-\alpha n\mathbf{1}_{[\tau _{nm}^{\ast }<T]}\mathbf{1}_{\overline{A}%
}+Y_{\tau _{nm}^{\ast }}^{u_{\tau _{nm}^{\ast }}}\mathbf{1}_{[\tau
_{nm}^{\ast }<T]}\right] .
\end{eqnarray*}%
Then the right-hand side converge to $-\infty $ as $n\rightarrow
\infty $. This contradicts the fact that $Y^{i}$ belong to
$\mathcal{S}^{2}$ and $\psi
_{i}(.,X)\in \mathcal{M}^{2,1}$. Henceforth the strategy is finite: $P(A)=0.$%

\subsection{Existence of processes $Y^{i}$, $i=1,\ldots ,m$}

The issue of existence of the processes $Y^{1},\ldots ,Y^{m}$ which satisfy (%
\ref{eqvt}) is also addressed in \cite{DHP}. Also for $n\geq 0$ let
us define the processes $(Y^{1,n},\ldots ,Y^{m,n})$ recursively as
follows: for $i\in \mathcal{I}$ we set,
\begin{equation}
Y_{t}^{i,0}=E\left[ \int_{t}^{T}\psi
_{i}(s,X_{s})ds|\mathcal{F}_{t}\right] ,\,\,0\leq t\leq T,
\label{y0}
\end{equation}%
and for $n\geq 1$,
\begin{equation}
Y_{t}^{i,n}=\limfunc{esssup}_{\tau \geq t}E\left[ \int_{t}^{\tau
}\psi _{i}(s,X_{s})ds+\max\limits_{k\in
\mathcal{I}^{-i}}(-g_{ik}(\tau ,X_{\tau })+Y_{\tau
}^{k,n-1})\mathbf{1}_{[\tau <T]}|\mathcal{F}_{t}\right] ,\,\,0\leq
t\leq T.  \label{eq24}
\end{equation}%
Then the sequence of processes $((Y^{1,n},\ldots ,Y^{m,n}))_{n\geq
0}$ have the following properties:

\begin{pro}
(\cite{DHP}, Pro. 3 and Th. 2)

\begin{itemize}
\item[$(i)$] for any $i\in \mathcal{I}$ and $n\geq 0$, the processes $%
Y^{1,n},\ldots ,Y^{m,n}$ are well defined, continuous and belong to $%
\mathcal{S}^{2}$, and verify
\begin{equation}
\forall \;t\leq T,\,\,Y_{t}^{i,n}\leq Y_{t}^{i,n+1}\leq E\left[
\int_{t}^{T}\left\{ \max_{1\leq i\leq m}|\psi _{i}(s,X_{s})|\right\} ds|%
\mathcal{F}_{t}\right] ;  \label{croi}
\end{equation}

\item[$(ii)$] there exist $m$ processes $Y^1,\ldots ,Y^m$ of $\mathcal{S}^2$
such that for any $i\in \mathcal{I}$:

\begin{itemize}
\item[$(a)$] $\forall \;t\leq T$, $Y_{t}^{i}=\lim_{n\rightarrow \infty
}\nearrow Y_{t}^{i,n}$ and
\begin{equation*}
E\left[ \sup_{s\leq T}|Y_{s}^{i,n}-Y_{s}^{i}|^{2}\right] \rightarrow 0\quad %
\mbox{ as }\quad n\rightarrow +\infty
\end{equation*}

\item[$(b)$] $\forall \;t\leq T$,
\begin{equation}
{Y}_{t}^{i}=\limfunc{esssup}_{\tau \geq t}E\left[ \int_{t}^{{\tau
}}\psi
_{i}(s,X_{s})ds+\max\limits_{k\in \mathcal{I}^{-i}}(-g_{ik}(\tau ,X_{\tau })+%
{Y}_{\tau }^{k})\mathbf{1}_{[\tau <T]}|\mathcal{F}_{t}\right]
\label{eq26}
\end{equation}%
i.e. ${Y}^{1},\ldots ,{Y}^{m}$ satisfy the Verification Theorem \ref%
{thmverif} ;

\item[$(c)$] $\forall \;t\leq T$,
\begin{equation}
{Y}_{t}^{i}=\limfunc{esssup}_{(\delta ,u)\in
\mathcal{D}_{t}^{i}}E\left[ \left\{ \int_{t}^{T}\psi
_{u_{s}}(s,X_{s})ds-\sum_{n\geq 1}g_{u_{\tau _{n-1}}u_{\tau
_{n}}}(\tau _{n},X_{\tau _{n}})\mathbf{1}_{[\tau _{n}<T]}\right\}
|\mathcal{F}_{t}\right]  \label{eq27}
\end{equation}%
where $\mathcal{D}_{t}^{i}=\{(\delta ,\xi )=((\tau _{n})_{n\geq
1},(\xi _{n})_{n\geq 1})\mbox { such that }u_{0}=i\mbox{ and }\tau
_{1}\geq t\}$. This characterization means that if at time $t$ the
production activity is in its regime $i$ then the optimal expected
profit is $Y_{t}^{i}$.

\item[$(d)$] the processes $Y^{1},\ldots ,Y^{m}$ verify the dynamical
programming principle of the $m$-states optimal switching problem, $i.e.$, $%
\forall \;t\leq T$,
\begin{equation}
Y_{t}^{i}=\limfunc{esssup}_{(\delta ,u)\in
\mathcal{D}_{t}^{i}}E\left[ \int_{t}^{{\tau _{n}}}\psi
_{u_{s}}(s,X_{s})ds-\sum_{1\leq k\leq
n}g_{u_{\tau _{k-1}}u_{{\tau _{k}}}}({\tau }_{k},X_{{\tau }_{k}})\mathbf{1}%
_{[\tau _{k}<T]}+\mathbf{1}_{[\tau _{n}<T]}Y_{\tau _{n}}^{u_{\tau _{n}}}|%
\mathcal{F}_{t}\right]  \label{dpp}
\end{equation}
\end{itemize}
\end{itemize}
\end{pro}

Note that except $(ii-d)$, the proofs of the other points are given in \cite%
{DHP}. The proof of $(ii.-d)$ can be easily deduced using relation (\ref%
{eq26}). From (\ref{eq26}) for any $i\in \mathcal{I}$, $t\in \lbrack
0,T]$ and $(\delta ,\xi )\in \mathcal{D}_{t}^{i}$ we have:
\begin{equation}
Y_{t}^{i}\geq E\left[ \int_{t}^{\tau _{n}}\psi
_{u_{s}}(s,X_{s})ds-\sum_{1\leq k\leq n}g_{u_{\tau _{k-1}}u_{u_{\tau _{k}}}}(%
{\tau }_{k},X_{{\tau }_{k}})\mathbf{1}_{[\tau
_{k}<T]}+\mathbf{1}_{[\tau _{n}<T]}Y_{\tau _{n}}^{u_{\tau
_{n}}}|\mathcal{F}_{t}\right] .  \label{eq28}
\end{equation}%
Next using the optimal strategy we obtain the equality instead of
inequality in (\ref{eq28}). Therefore the relation (\ref{dpp}) holds
true. $\Box $

\begin{rem}
\label{unic}Note that the characterization (\ref{eq27}) implies that
the processes $Y^1,\ldots ,Y^m$ of $\mathcal{S}^2$ which satisfy the
Verification Theorem are unique.
\end{rem}

\section{Existence of a solution for the system of variational inequalities}

\subsection{Connection with BSDEs with one reflecting barrier}

Let $(t,x)\in [0,T]\times \mathbb{R}^k$ and let $(X^{t,x}_s)_{s\leq
T}$ be the solution of the following standard SDE:
\begin{equation}  \label{sde}
dX^{t,x}_s=b(s,X_s^{t,x})ds+\sigma(s,X_s^{t,x})dB_s \mbox{ for }t\leq s\leq T%
\mbox{ and }X_s^{t,x}=x \mbox{ for }s\leq t
\end{equation}%
where the functions $b$ and $\sigma$ are the ones of (\ref{regbs1}).
These
properties of $\sigma$ and $b$ imply in particular that the process $%
(X^{t,x}_s)_{0\le s\leq T}$ solution of the standard SDE (\ref{sde})
exists and is unique, for any $t\in [0, T]$ and $x\in \mathbb{R}^k$.

The operator $\mathcal{A}$ that is appearing in (\ref{generateur})
is the infinitesimal generator associated with $X^{t,x}$. In the
following result we collect some properties of $X^{t,x}$.

\begin{pro}
\label{estimx} (\cite{RY}) The process $X^{t,x}$ satisfies the
following estimates:

\begin{itemize}
\item[$(i)$] For any $q\geq 2$, there exists a constant $C$ such that
\begin{equation}
E\left[ \sup_{0\leq s\leq T}|X_{s}^{t,x}|^{q}\right] \leq
C(1+|x|^{q}). \label{estimat1}
\end{equation}

\item[$(ii)$] There exists a constant $C$ such that for any $t,t^{\prime
}\in \lbrack 0,T]$ and $x,x^{\prime }\in \mathbb{R}^{k}$,
\begin{equation}
E\left[ \sup_{0\leq s\leq T}|X_{s}^{t,x}-X_{s}^{t^{\prime }x^{\prime }}|^{2}%
\right] \leq C(1+|x|^{2})(\left\vert x-x^{\prime }\right\vert
^{2}+|t-t^{\prime }|).\Box  \label{estimat2}
\end{equation}
\end{itemize}
\end{pro}

We consider a BSDE with one reflecting barrier introduced in
\cite{EKal}. This notion will allow us to make the connection
between the variational inequalities (\ref{sysvi1}) and the
$m$-states optimal switching problem described in the previous
section.

\medskip

Let $f:[0,T]\times \mathbb{R}^{k+1+d}\rightarrow \mathbb{R}$,
$h:[0,T]\times \mathbb{R}^{k}\rightarrow \mathbb{R}$ and
$g:\mathbb{R}^{k}\rightarrow \mathbb{R}$ be continuous, of
polynomial growth and such that $h(x,T)\leq
g(x)$. Moreover we assume that for any $(t,x)\in \lbrack 0,T]\times \mathbb{R%
}^{k}$, the mapping $(y,z)\in \mathbb{R}^{1+d}\mapsto f(t,x,y,z)$ is
uniformly Lipschitz. Then we have the following result related to
BSDEs with one reflecting barrier:

\begin{theo}
(\cite{EKal}, Th. 5.2 and 8.5) \label{theo-elk}For any $(t,x)\in \lbrack 0,T]\times \mathbb{R%
}^{k}$, there exits a unique triple of processes
$(Y^{t,x},Z^{t,x},K^{t,x})$ such that:
\begin{equation}
\left\{
\begin{array}{l}
Y^{t,x},K^{t,x}\in \mathcal{S}^{2}\mbox{ and }Z^{t,x}\in \mathcal{M}%
^{2,d};\,K^{t,x}\mbox{ is  non-decreasing
and }K_{0}^{t,x}=0, \\
\\
Y_{s}^{t,x}=g(X_{T}^{t,x})+\dint%
\nolimits_{s}^{T}f(r,X_{r}^{t,x},Y_{r}^{t,x},Z_{r}^{t,x})dr-\dint%
\nolimits_{s}^{T}Z_{r}^{t,x}dB_{r}+K_{T}^{t,x}-K_{s}^{t,x},\,\,s\leq T \\
\\
Y_{s}^{t,x}\geq h(s,X_{s}^{t,x}),\,\forall \;s\leq T\mbox{ and }%
\dint\nolimits_{0}^{T}(Y_{r}^{t,x}-h(r,X_{r}^{t,x}))dK_{r}^{t,x}=0.%
\end{array}%
\right.
\end{equation}%
Moreover, the following characterization of $Y^{t,x}$ as a Snell
envelope holds true:
\begin{equation}
\forall \;s\leq T,\,\,Y_{s}^{t,x}=\limfunc{esssup}_{\tau \in \mathcal{T}%
_{t}}E\left[ \int_{t}^{\tau
}f(r,X_{r}^{t,x},Y_{r}^{t,x},Z_{r}^{t,x})dr+h(\tau ,X_{\tau }^{t,x})\mathbf{1%
}_{[\tau <T]}+g(X_{T}^{t,x})\mathbf{1}_{[\tau
=T]}|\mathcal{F}_{s}\right] . \label{snellenv}
\end{equation}%
There exists a deterministic continuous function $u:[0,T]\times \mathbb{R}%
^{k}\rightarrow \mathbb{R}$ with polynomial growth such that:%
\begin{equation*}
\forall \;s\in \lbrack t,T],Y_{s}^{t,x}=u(s,X_{s}^{t,x}).
\end{equation*}%
and the function $u$ is the unique viscosity solution in the class
of continuous function with polynomial growth of the following PDE
with
obstacle:%
\begin{equation*}
\left\{
\begin{array}{l}
\min \{u(t,x)-h(t,x),-\partial _{t}u(t,x)-\mathcal{A}u(t,x)-f(t,x,u(t,x),%
\sigma (t,x)^{\ast }\nabla u(t,x))\}=0, \\
u(T,x)=g(x).\Box%
\end{array}%
\right.
\end{equation*}
\end{theo}

\subsection{Existence of a solution for the system of variational
inequalities}

Let $(Y^{1,t,x}_s,\ldots ,Y^{m,t,x}_s)_{0\le s\leq T}$ be the
processes which satisfy the Verification Theorem \ref{thmverif} in
the case when the
process $X\equiv X^{t,x}$. Therefore using the characterization (\ref%
{snellenv}), there exist processes $K^{i,t,x}$

and $Z^{i,t,x}$, $i\in \mathcal{I}$, such that the triples ($%
Y^{i,t,x},Z^{i,t,x},K^{i,t,x})$ are unique solutions of the
following reflected BSDEs: for any $i=1,\ldots ,m$ we have
\begin{equation}
\left\{
\begin{array}{l}
Y^{i,t,x},K^{i,t,x}\in \mathcal{S}^{2}\mbox{ and }Z^{i,t,x}\in \mathcal{M}%
^{2,d};\,K^{i,t,x}\mbox{ is
non-decreasing and }K_{0}^{i,t,x}=0, \\
\\
Y_{s}^{i,t,x}=\dint\nolimits_{s}^{T}\psi
_{i}(r,X_{r}^{t,x})du-\dint%
\nolimits_{s}^{T}Z_{r}^{i,t,x}dB_{r}+K_{T}^{i,t,x}-K_{s}^{i,t,x},\,\,\,0\leq
s\leq T,\,\,Y_{T}^{i,t,x}=0, \\
\\
Y_{s}^{i,t,x}\geq \max\limits_{j\in \mathcal{I}%
^{-i}}(-g_{ij}(s,X_{s}^{t,x})+Y_{s}^{j,t,x}),\,\,0\leq s\leq T, \\
\\
\dint\nolimits_{0}^{T}\left( Y_{r}^{i,t,x}-\max\limits_{j\in \mathcal{I}%
^{-i}}(-g_{ij}(r,X_{r}^{t,x})+Y_{r}^{j,t,x})\right) dK_{r}^{i,t,x}=0.%
\end{array}%
\right.
\end{equation}%
Moreover we have the following representation of $Y$.

\begin{pro}
There are deterministic functions $v^{1},\ldots ,v^{m}$
$:[0,T]\times \mathbb{R}^{k}\rightarrow \mathbb{R}$ such that:
\begin{equation*}
\forall \;(t,x)\in \lbrack 0,T]\times \mathbb{R}^{k},\forall \;s\in
\lbrack t,T],Y_{s}^{i,t,x}=v^{i}(s,X_{s}^{t,x}),\,\,i=1,\ldots ,m.
\end{equation*}%
and the functions $v^{i}$, $i=1,\ldots ,m,$ are lower
semi-continuous and of polynomial growth.
\end{pro}

$Proof$:

For $n\geq 0$ let $(Y_{s}^{n,1,t,x},\ldots ,Y_{s}^{n,m,t,x})_{0\leq
s\leq T}$ be the processes constructed in (\ref{y0})-(\ref{eq24}).
Therefore using an induction argument and Theorem \ref{theo-elk}
there exist deterministic continuous with polynomial growth
functions $v^{i,n}$ ($i=1,\ldots ,m$) such that for any $(t,x)\in
\lbrack 0,T]\times \mathbb{R}^{k}$, $\forall \;s\in
\lbrack t,T]$, $Y_{s}^{n,i,t,x}=v^{i,n}(s,X_{s}^{t,x})$. Inequality (\ref%
{croi}) yields
\begin{equation*}
Y_{t}^{n,i,t,x}\leq Y_{t}^{n+1,i,t,x}\leq E\left[
\int_{t}^{T}\left\{
\max_{1\leq i\leq m}|\psi _{i}(s,X_{s}^{t,x})|\right\} ds|\mathcal{F}_{t}%
\right]
\end{equation*}%
since $Y_{t}^{n,i,t,x}$ is deterministic. Therefore combining the
polynomial growth of $\psi _{i}$ and estimate (\ref{estimat1}) for
$X^{t,x}$ we obtain:
\begin{equation*}
v^{i,n}(t,x)\leq v^{i,n+1}(t,x)\leq C(1+|x|^{p})
\end{equation*}%
for some constants $C$ and $p$ independent of $n$. In order to
complete the proof it is enough to set
$v^{i}(t,x):=\lim_{n\rightarrow \infty
}v^{i,n}(t,x),(t,x)\in \lbrack 0,T]\times \mathbb{R}^{k}$ since $%
Y^{i,n,t,x}\nearrow Y^{i,t,x}$ as $n\rightarrow \infty $. $\Box $

We are now going to focus on the continuity of the functions
$v^1,...,v^m$. But first let us deal with some properties of the
optimal strategy which exist thanks to Theorem \ref{thmverif}.

\begin{pro}
\label{optimal-s} Let $(\delta ,u)=((\tau _{n})_{n\geq 1},(\xi
_{n})_{n\geq
1})$ be an optimal strategy finite, then there exist two positive constant $%
C $ and $p$ which do not depend on $t$ and $x$ such that:

\begin{itemize}
\item[$(i)$] if $E\left[ \sum\limits_{k\geq 1}g_{u_{\tau _{k}^{\ast
}}u_{\tau _{k-1}^{\ast }}}({\tau _{k}^{\ast }},X_{{\tau _{k}^{\ast }}}^{t,x})%
\mathbf{1}_{[\tau _{k}^{\ast }<T]}\right] <+\infty $, then
\begin{equation}
\forall \;n\geq 1,\,\,P[\tau _{n}<T]\leq
\frac{C(1+|x|^{p})}{\sqrt{n}}. \label{estiopt1}
\end{equation}

\item[$(ii)$] If $g_{ij}$ is constant, then
\begin{equation}
\forall \;n\geq 1,\,\,P[\tau _{n}<T]\leq \frac{mC(1+|x|^{p})}{n}.
\label{estiopt2}
\end{equation}
\end{itemize}
\end{pro}

$Proof.$ (i) We will show by contradiction, suppose $\exists \
n_{0},\,\forall \;n_{1}\geq n_{0},\,\,P[\tau _{n_{1}}<T]>\frac{C(1+|x|^{p})}{%
\sqrt{n_{1}}}.$

Recall the characterization of (\ref{eq27}) that reads as:
\begin{equation*}
Y_{0}^{1,t,x}=\limfunc{esssup}_{(\delta ,u)\in \mathcal{D}}E\left[
\int_{0}^{T}\psi _{u_{r}}(r,X_{r}^{t,x})dr-\sum_{k\geq 1}g_{u_{\tau
_{k-1}}u_{\tau _{k}}}(\tau _{k},X_{\tau
_{k}}^{t,x})\mathbf{1}_{[\tau _{k}<T]}\right] .
\end{equation*}

If $(\delta ,u)=((\tau _{n})_{n\geq 1},(\xi _{n})_{n\geq 1})$ is the
optimal strategy then
\begin{eqnarray*}
&&Y_{0}^{1,t,x}-E\left[ \sum\limits_{k\geq 1}g_{u_{\tau _{k}}u_{\tau
_{k-1}}}({\tau _{k}},X_{{\tau _{k}}}^{t,x})\mathbf{1}_{[\tau
_{k}<T]}\right]
\\
&& \\
&=&E\left[ \int_{0}^{T}\psi _{u_{r}}(r,X_{r}^{t,x})dr-\sum_{k\geq
1}g_{u_{\tau _{k-1}}u_{\tau _{k}}}(\tau _{k},X_{\tau _{k}}^{t,x})\mathbf{1}%
_{[\tau _{k}<T]}\right] \\
&& \\
&&-E\left[ \sum\limits_{k\geq 1}g_{u_{\tau _{k}}u_{\tau _{k-1}}}({\tau _{k}}%
,X_{{\tau _{k}}}^{t,x})\mathbf{1}_{[\tau _{k}<T]}\right] .
\end{eqnarray*}%
Taking into account that $g_{ij}+g_{ji}>\alpha $ for any $i\neq j$
and for any $k\leq n_{1}$, $[\tau _{n_{1}}<T]\subset \lbrack \tau
_{k}<T]$ we obtain:
\begin{eqnarray*}
&&Y_{0}^{1,t,x}-E\left[ \sum\limits_{k\geq 1}g_{u_{\tau _{k}}u_{\tau
_{k-1}}}({\tau _{k}},X_{{\tau _{k}}}^{t,x})\mathbf{1}_{[\tau
_{k}<T]}\right]
\\
&& \\
&\leq &E\left[ \int_{0}^{T}\psi _{u_{r}}(r,X_{r}^{t,x})dr\right]
-n_{1}\alpha P[\tau _{n_{1}}<T] \\
&& \\
&\leq &E\left[ \int_{0}^{T}\psi _{u_{r}}(r,X_{r}^{t,x})dr\right]
-n_{1}\alpha \frac{C(1+|x|^{p})}{\sqrt{n_{1}}}.
\end{eqnarray*}%
As $n_{1}$ is arbitrary then putting $n_{1}\rightarrow +\infty $ to
obtain:
\begin{equation*}
Y_{0}^{1,t,x}-E\left[ \sum\limits_{k\geq 1}g_{u_{\tau _{k}}u_{\tau _{k-1}}}({%
\tau _{k}},X_{{\tau _{k}}}^{t,x})\mathbf{1}_{[\tau _{k}<T]}\right]
\leq -\infty ,
\end{equation*}%
which is a contradiction.\newline (ii) If $(\delta ,u)=((\tau
_{n})_{n\geq 1},(\xi _{n})_{n\geq 1})$ is the optimal strategy and
$g_{ij}$ is constant then we have:
\begin{eqnarray*}
Y_{0}^{1,t,x} &=&E\left[ \int_{0}^{T}\psi
_{u_{r}}(r,X_{r}^{t,x})dr-\sum_{k\geq 1}g_{u_{\tau _{k-1}}u_{\tau _{k}}}%
\mathbf{1}_{[\tau _{k}<T]}\right] \\
&& \\
&\leq &E\left[ \int_{0}^{T}\psi
_{u_{r}}(r,X_{r}^{t,x})dr-\sum_{k=1}^{nm}g_{u_{\tau _{k-1}}u_{\tau _{k}}}%
\mathbf{1}_{[\tau _{k}<T]}\right] .
\end{eqnarray*}%
From (\ref{somme}) we have:
\begin{equation*}
Y_{0}^{1,t,x}\leq E\left[ \int_{0}^{T}\psi
_{u_{r}}(r,X_{r}^{t,x})dr\right] -\alpha nP[\tau _{nm}<T].
\end{equation*}%
Then,
\begin{eqnarray*}
n\alpha P[\tau _{nm} &<&T]\leq E\left[ \int_{0}^{T}\mid \psi
_{u_{r}}(r,X_{r}^{t,x})\mid dr\right] -Y_{0}^{1,t,x} \\
&& \\
&\leq &E\left[ \int_{0}^{T}\mid \psi _{u_{r}}(r,X_{r}^{t,x})\mid
dr\right] -Y_{0}^{0,1,t,x}.
\end{eqnarray*}

\begin{rem}
The estimate (\ref{estiopt1}) and (\ref{estiopt2}) are also valid
for the optimal strategy if at the initial time the state of the
plant is an arbitrary $i\in \mathcal{I}$. $\Box $
\end{rem}

\medskip

Next,for $i\in \mathcal{I}$, let $%
(y_{s}^{i,t,x},z_{s}^{i,t,x},k_{s}^{i,t,x})_{0\leq s\leq T}$ be the
processes defined as follows:
\begin{equation}
\left\{
\begin{array}{l}
y^{i,t,x},k^{i,t,x}\in \mathcal{S}^{2}\mbox{ and }z^{i,t,x}\in \mathcal{M}%
^{2,d};\,k^{i,t,x}\mbox{ is  non-decreasing and }k_{0}^{i,t,x}=0, \\
\\
y_{s}^{i,t,x}=\dint\nolimits_{s}^{T}\psi _{i}(r,X_{r}^{t,x})\mathbf{1}%
_{[r\geq
t]}dr-\dint\nolimits_{s}^{T}z_{r}^{i,t,x}dB_{r}+k_{T}^{i,t,x}-k_{s}^{i,t,x},%
\,\,\,0\leq s\leq T,\,\,y_{T}^{i,t,x}=0 \\
\\
y_{s}^{i,t,x}\geq l_{s}^{i,t,x}:=\max\limits_{j\in \mathcal{I}%
^{-i}}\{-g_{ij}(t\vee s,X_{t\vee
s}^{t,x})+y_{s}^{j,t,x})\},\,\,\,\forall
\;s\leq T, \\
\\
\dint\nolimits_{0}^{T}(y_{r}^{i,t,x}-l_{r}^{i,t,x})dk_{r}^{i,t,x}=0.%
\end{array}%
\right.  \label{y2}
\end{equation}%
The existence of $(y^{i,t,x},z^{i,t,x},k^{i,t,x}),i\in \mathcal{I}$,
is obtained in the same way as the one of
$(Y^{i,t,x},Z^{i,t,x},K^{i,t,x})$. By uniqueness we obtain for any
$(t,x)\in \lbrack 0,T]\times \mathbb{R}^{k}$,
for any $s\in \lbrack 0,t]$ we have $y_{s}^{i,t,x}=Y_{t}^{i,t,x}$, $%
z_{s}^{i,t,x}=0$ and $k_{s}^{i,t,x}=0$. \medskip

We are now ready to give the continuity of the value functions, when
the strategy optimal is finite.

\begin{theo}
The functions $(v^1,\ldots ,v^m):[0,T]\times \mathbb{R}^k\rightarrow \mathbb{%
R}$ are continuous and solution in viscosity sense of the system of
variational inequalities with inter-connected obstacles
(\ref{sysvi1}).
\end{theo}

$Proof.$ The continuity of the value functions follows from the
dynamic programming principle and is proved in \cite{EH}.$\Box$

\section{Uniqueness of the solution of the system}

In this section we address the main question of this paper, that is
uniqueness of the viscosity solution of the system (\ref{sysvi1}).

\begin{theo}
\label{uni}The solution in viscosity sense of the system of
variational inequalities with inter-connected obstacles
(\ref{sysvi1}) is unique in the space of continuous functions on
$[0,T]\times \mathbb{R}^{k}$ which satisfy a polynomial growth
condition, i.e. in the space
\begin{equation*}
\mathcal{C}:=\left\{ \varphi :[0,T]\times \mathbb{R}^{k}\rightarrow \mathbb{R%
},\mbox{
continuous and for any }(t,x),\,|\varphi (t,x)|\leq C(1+|x|^{\gamma })%
\mbox{ for some constants }C\mbox{ and }\gamma \right\} .
\end{equation*}
\end{theo}
$Proof.$ The proof is divided in four steps. We will show by
contradiction that if $u_{1},\ldots ,u_{m}$ and $w_{1},\ldots
,w_{m}$ are a subsolution and a supersolution respectively for
(\ref{sysvi1}) then for any $i=1,\ldots ,m$, $u_{i}\leq w_{i}$.
Therefore if we have two solutions of (\ref{sysvi1})
then they are obviously equal. Actually for some $R>0$ suppose there exists $%
(\overline{t},\overline{x},\overline{i})\in (0,T)\times B_{R}\times \mathcal{%
I}$ $(B_{R}:=\{x\in \mathbb{R}^{k};|x|<R\})$ such that:
\begin{equation}
\max\limits_{(t,x,i)}(u_{i}(t,x)-w_i(t,x))=u_{\overline{i}}(\overline{t},%
\overline{x})-w_{\overline{i}}(\overline{t},\overline{x})=\eta >0.
\label{comp-equ}
\end{equation}%
\textbf{Step 1. } Let us take $\theta ,\lambda $ and $\beta \in
(0,1]$ small enough, so that the following holds:
\begin{equation}
\left\{
\begin{array}{l}
\beta T<\dfrac{\eta }{6} \\
2\theta |\overline{x}|^{2\gamma +2}<\dfrac{\eta }{6} \\
-\lambda w_{\overline{i}}(\overline{t},\overline{x})<\dfrac{\eta }{6} \\
\dfrac{\lambda }{\overline{t}}<\dfrac{\eta }{6} \\
\lambda \max\limits_{i\neq j}g_{ji}(\overline{t},\overline{x})<\dfrac{\eta }{%
6}.%
\end{array}%
\right.
\end{equation}%
Here $\gamma $ is the growth exponent of the functions which w.l.o.g
we assume integer and $\geq 2$. Then, for a small $\varepsilon >0$,
let us define:
\begin{equation}
\Phi _{\varepsilon }^{i}(t,x,y)=u_{i}(t,x)-w_{i}^{\lambda }(t,y)-\frac{1}{%
2\varepsilon }(|x-y|^{2\gamma }+|x-\overline{x}|^{2\gamma }+(t-\overline{t}%
)^{2\gamma })-\theta (|x|^{2\gamma +2}+|y|^{2\gamma +2})+\beta t-\frac{%
\lambda }{t}  \label{phi}
\end{equation}%
where, $w_{i}^{\lambda }(t,x)=(1-\lambda )w(t,x)+\lambda \alpha _{i}$ and $%
\alpha _{i}=\min_{j\in
\mathcal{I}^{-i}}g_{ji}(\overline{t},\overline{x})$.
By the growth assumption on $u_{i}$ and $w_{i}$, there exists a $%
(t_{0},x_{0},y_{0},i_{0})\in (0,T]\times B_{R}\times B_{R}\times \mathcal{I}$%
, such that:
\begin{equation*}
\Phi _{\varepsilon
}^{i_{0}}(t_{0},x_{0},y_{0})=\max\limits_{(t,x,y,i)}\Phi
_{\varepsilon }^{i}(t,x,y).
\end{equation*}%
On the other hand, from $\Phi _{\varepsilon
}^{i_{0}}(t_{0},x_{0},y_{0})\geq \Phi _{\varepsilon
}^{i_{0}}(\overline{t},\overline{x},\overline{x})$, we have
\begin{eqnarray}
\frac{1}{2\varepsilon }(|x_{0}-y_{0}|^{2\gamma }+|x_{0}-\overline{x}%
|^{2\gamma }+(t_{0}-\overline{t})^{2\gamma }) &\leq
&(u_{i_{0}}(t_{0},x_{0})-u_{i_{0}}(\overline{t},\overline{x}%
))+(w_{i_{0}}^{\lambda
}(\overline{t},\overline{x})-w_{i_{0}}^{\lambda
}(t_{0},y_{0} \\
&&-\theta (|x_{0}|^{2\gamma +2}+|y_{0}|^{2\gamma +2}-2|\overline{x}%
|^{2\gamma +2})+\beta (t_{0}-\overline{t})-\frac{\lambda }{t_{0}}+\frac{%
\lambda }{\overline{t}}  \notag
\end{eqnarray}%
and consequently $\frac{1}{2\varepsilon }(|x_{0}-y_{0}|^{2\gamma }+|x_{0}-%
\overline{x}|^{2\gamma }+(t_{0}-\overline{t})^{2\gamma })$ is
bounded, and
as $\varepsilon \rightarrow 0$, $|x_{0}-y_{0}|\rightarrow 0$, $|x_{0}-%
\overline{x}|\rightarrow 0$ and $(t_{0}-\overline{t})\rightarrow 0$. Since $%
u_{i_{0}}$ and $w_{i_{0}}^{\lambda }$ are uniformly continuous on $%
[0,T]\times \overline{B}_{R}$, then $\frac{1}{2\varepsilon }%
(|x_{0}-y_{0}|^{2\gamma }+|x_{0}-\overline{x}|^{2\gamma }+(t_{0}-\overline{t}%
)^{2\gamma })$ as $\varepsilon \rightarrow 0.$\\
\textbf{Step 2. } We now show that $t_{0}<T.$ If $t_{0}=T$ then,
\begin{equation*}
\Phi _{\varepsilon }^{\overline{i}}(\overline{t},\overline{x},\overline{x}%
)\leq \Phi _{\varepsilon }^{i_{0}}(T,x_{0},y_{0}),
\end{equation*}%
and,
\begin{equation*}
u_{\overline{i}}(\overline{t},\overline{x})-(1-\lambda )w_{\overline{i}}(%
\overline{t},\overline{x})-2\theta |\overline{x}|^{2\gamma +2}+\beta
\overline{t}-\frac{\lambda }{\overline{t}}\leq \lambda \alpha
_{i_{0}}+\beta T-\frac{\lambda }{T},
\end{equation*}%
since $u_{i_{0}}(T,x_{0})=w_{i_{0}}(T,y_{0})=0$ and $\alpha _{\overline{i}%
}\geq 0$. Then thanks to (\ref{comp-equ}) we have,
\begin{eqnarray*}
\eta &\leq &-\lambda
w_{\overline{i}}(\overline{t},\overline{x})+\lambda
\alpha _{i_{0}}+\beta T+2\theta |\overline{x}|^{2\gamma +2}+\frac{\lambda }{%
\overline{t}} \\
\eta &<&\frac{5}{6}\eta .
\end{eqnarray*}%
which yields a contradiction and we have $t_{0}\in (0,T)$. \newline
\textbf{Step 3. }We now claim that:
\begin{equation}
u_{i_{0}}(t_{0},x_{0})-\max\limits_{j\in \mathcal{I}^{-i_{0}}}%
\{-g_{i_{0}j}(t_{0},x_{0})+u_{j}(t_{0},x_{0})\}>0.
\label{visco-comp1}
\end{equation}%
Indeed if
\begin{equation*}
u_{i_{0}}(t_{0},x_{0})-\max\limits_{j\in \mathcal{I}^{-i_{0}}}%
\{-g_{i_{0}j}(t_{0},x_{0})+u_{j}(t_{0},x_{0})\}\leq 0
\end{equation*}%
then there exists $k\in \mathcal{I}^{-i_{0}}$ such that:
\begin{equation*}
u_{i_{0}}(t_{0},x_{0})\leq
-g_{i_{0}k}(t_{0},x_{0})+u_{k}(t_{0},x_{0}).
\end{equation*}%
Now, we then see that
\begin{equation}
\begin{array}{l}
w_{i_{0}}^{\lambda }(t_{0},y_{0})-\max\limits_{j\in \mathcal{\ I}%
^{-i_{0}}}(-g_{i_{0}j}(t_{0},y_{0})+w_{j}^{\lambda }(t_{0},y_{0})) \\
\quad =\lambda \alpha _{i_{0}}+(1-\lambda
)w_{i_{0}}(t_{0},y_{0})-\max\limits_{j\in
\mathcal{I}^{-i_{0}}}[(1-\lambda
)(-g_{i_{0}j}(t_{0},y_{0})+w_{j}(t_{0},y_{0}))+\lambda \alpha
_{j}-\lambda
g_{i_{0}j}(t_{0},y_{0})] \\
\quad \geq (1-\lambda )[w_{i_{0}}(t_{0},y_{0})-\max\limits_{j\in \mathcal{I}%
^{-i_{0}}}(-g_{i_{0}j}(t_{0},y_{0})+w_{j}(t_{0},y_{0}))]+\lambda
(\alpha _{i_{0}}-\max\limits_{j\in \mathcal{I}^{-i_{0}}}(\alpha
_{j}-g_{i_{0}j}(t_{0},y_{0})) \\
\quad \geq \lambda (\alpha _{i_{0}}+\min\limits_{j\in \mathcal{I}%
^{-i_{0}}}(g_{i_{0}j}(t_{0},y_{0})-\alpha _{j}),%
\end{array}
\label{super-ine}
\end{equation}%
let $i_{1}\in \mathcal{I}^{-i_{0}}$ such that $\alpha
_{i_{0}}=g_{i_{1}i_{0}}(\overline{t},\overline{x})$ and set
$i_{2}\in \mathcal{I}^{-i_{0}}$ such that
\begin{equation*}
\min\limits_{j\in
\mathcal{I}^{-i_{0}}}(g_{i_{0}j}(t_{0},y_{0})-\alpha
_{j})=g_{i_{0}i_{2}}(t_{0},y_{0})-\alpha _{i_{2}}.
\end{equation*}%
Then we have
\begin{equation*}
\begin{array}{l}
\alpha _{i_{0}}+\min\limits_{j\in \mathcal{I}%
^{-i_{0}}}(g_{i_{0}j}(t_{0},y_{0})-\alpha _{j})=g_{i_{1}i_{0}}(\overline{t},%
\overline{x})+g_{i_{0}i_{2}}(t_{0},y_{0})-\alpha _{i_{2}} \\
=g_{i_{1}i_{0}}(\overline{t},\overline{x})+g_{i_{0}i_{2}}(\overline{t},%
\overline{x})-g_{i_{0}i_{2}}(\overline{t},\overline{x}%
)+g_{i_{0}i_{2}}(t_{0},y_{0})-\alpha _{i_{2}} \\
=\nu -g_{i_{0}i_{2}}(\overline{t},\overline{x})+g_{i_{0}i_{2}}(t_{0},y_{0})%
\end{array}%
\end{equation*}

where
\begin{equation*}
\nu=g_{i_1i_0}(\overline{t},\overline{x})+g_{i_0i_2}(\overline{t},\overline{x%
})-\min\limits_{j\in\mathcal{I}^{-i_2}}g_{ji_2}(\overline{t},\overline{x}%
)>g_{i_1i_2}(\overline{t},\overline{x})-\min\limits_{j\in\mathcal{I}%
^{-i_2}}g_{ji_2}(\overline{t},\overline{x})\geq 0.
\end{equation*}
From (\ref{super-ine}) we have
\begin{equation*}
w^\lambda_{i_0}(t_0,y_0)-(-g_{i_0
k}(t_0,y_0)+w^\lambda_k(t_0,y_0))\geq\lambda\nu -\lambda g_{i_0i_2}(\overline{t},%
\overline{x})+\lambda g_{i_0i_2}(t_0,y_0).
\end{equation*}
It follows that:
\begin{equation*}
u_{i_0}(t_0,x_0)- w^\lambda_{i_0}(t_0,y_0)
-(u_{k}(t_0,x_0)-w^\lambda_{k}(t_0,y_0))\leq -\lambda\nu +\lambda g_{i_0i_2}(%
\overline{t},\overline{x})-\lambda g_{i_0i_2}(t_0,y_0)+g_{i_0
k}(t_0,y_0)-g_{i_0 k}(t_0,x_0).
\end{equation*}
Now taking into account of (\ref{phi}) to obtain:
\begin{equation*}
\Phi^{i_0}_{\varepsilon}(t_{0},x_{0},y_{0})-\Phi^{k}_{%
\varepsilon}(t_{0},x_{0},y_{0})\leq -\lambda\nu +\lambda g_{i_0i_2}(\overline{t},%
\overline{x})-\lambda g_{i_0i_2}(t_0,y_0)+g_{i_0 k}(t_0,y_0)-g_{i_0
k}(t_0,x_0).
\end{equation*}
But this contradicts the definition of $i_0$, since $g_{i_0i_2}$ and
$g_{i_0
k}$ are uniformly continuous on $[0,T]\times \overline{B}_R$ and the claim (%
\ref{visco-comp1}) holds.\\ \textbf{Step 4. } To complete the proof
it remains to show contradiction. Let us denote
\begin{equation}
\varphi _{\varepsilon }(t,x,y)=\frac{1}{2\varepsilon }(|x-y|^{2\gamma }+|x-%
\overline{x}|^{2\gamma }+(t-\overline{t})^{2\gamma })+\theta
(|x|^{2\gamma +2}+|y|^{2\gamma +2})-\beta t+\frac{\lambda }{t}.
\end{equation}%
Then we have:
\begin{equation}
\left\{
\begin{array}{l}
\label{derive}D_{t}\varphi _{\varepsilon }(t,x,y)=-\beta -\dfrac{\lambda }{%
t^{2}}+\dfrac{\gamma }{\varepsilon }(t-\overline{t})(t-\overline{t}%
)^{2\gamma -2}, \\
\\
D_{x}\varphi _{\varepsilon }(t,x,y)=\dfrac{\gamma }{\varepsilon }%
(x-y)|x-y|^{2\gamma -2}+\dfrac{\gamma }{\varepsilon }(x-\overline{x})|x-%
\overline{x}|^{2\gamma -2}+\theta (2\gamma +2)x|x|^{2\gamma }, \\
\\
D_{y}\varphi _{\varepsilon }(t,x,y)=-\dfrac{\gamma }{\varepsilon }%
(x-y)|x-y|^{2\gamma -2}+\theta (2\gamma +2)y|y|^{2\gamma }, \\
\\
B(t,x,y)=D_{x,y}^{2}\varphi _{\varepsilon }(t,x,y)=\dfrac{1}{\varepsilon }%
\begin{pmatrix}
a_{1}(x,y) & -a_{1}(x,y) \\
-a_{1}(x,y) & a_{1}(x,y)%
\end{pmatrix}%
+%
\begin{pmatrix}
a_{2}(x)+a_{3}(x) & 0 \\
0 & a_{2}(y)%
\end{pmatrix}
\\
\\
\mbox{ with }a_{1}(x,y)=\gamma |x-y|^{2\gamma -2}I+\gamma (2\gamma
-2)(x-y)(x-y)^{\ast }|x-y|^{2\gamma -4}, \\
\\
a_{2}(x)=\theta (2\gamma +2)|x|^{2\gamma }I+2\theta \gamma (2\gamma
+2)xx^{\ast }|x|^{2\gamma -2}\mbox{ and
} \\
\\
a_{3}(x)=\dfrac{\gamma }{\varepsilon }|x-\overline{x}|^{2\gamma -2}I+\dfrac{%
\gamma (2\gamma -2)}{\varepsilon }(x-\overline{x})(x-\overline{x})^{\ast }|x-%
\overline{x}|^{2\gamma -4}.%
\end{array}%
\right.
\end{equation}%
Taking into account (\ref{visco-comp1}) then applying the result by
Crandall et al. (Theorem 8.3, {\cite{CIL}) to the function
\begin{equation*}
u_{i_{0}}(t,x)-w_{i_{0}}^{\lambda }(t,y)-\varphi _{\varepsilon
}(t,x,y)
\end{equation*}%
at the point $(t_{0},x_{0},y_{0})$, for any $\varepsilon _{1}>0$, we
can find $c,d\in \mathbb{R}$ and $X,Y\in \mathcal{S}_{k}$, such
that: }

\begin{equation}
\left\{
\begin{array}{l}
\left( c,\dfrac{\gamma }{\varepsilon
}(x_{0}-y_{0})|x_{0}-y_{0}|^{2\gamma
-2}+\dfrac{\gamma }{\varepsilon }(x_{0}-\overline{x})|x_{0}-\overline{x}%
|^{2\gamma -2}+\theta (2\gamma +2)x_{0}|x_{0}|^{2\gamma },X\right)
\in
J^{2,+}(u_{i_{0}}(t_{0},x_{0})), \\
\\
(-d,\dfrac{\gamma }{\varepsilon }(x_{0}-y_{0})|x_{0}-y_{0}|^{2\gamma
-2}-\theta (2\gamma +2)y_{0}|y_{0}|^{2\gamma },Y)\in
J^{2,-}(w_{i_{0}}^{\lambda }(t_{0},y_{0})), \\
\\
c+d=D_{t}\varphi _{\varepsilon }(t_{0},x_{0},y_{0})=-\beta -\dfrac{\lambda }{%
t_{0}^{2}}+\dfrac{\gamma }{\varepsilon }(t_{0}-\overline{t})(t_{0}-\overline{%
t})^{2\gamma -2}\mbox{ and finally } \\
\\
-\left( \dfrac{1}{\varepsilon _{1}}+||B(t_{0},x_{0},y_{0})||\right)
I\leq
\begin{pmatrix}
X & 0 \\
0 & -Y%
\end{pmatrix}%
\leq B(t_{0},x_{0},y_{0})+\varepsilon _{1}B(t_{0},x_{0},y_{0})^{2}.%
\end{array}%
\right.  \label{lemmeishii}
\end{equation}%
By (\ref{visco-comp1}), and the definition of viscosity solution, we
get:
\begin{eqnarray*}
-c-\frac{1}{2}Tr[\sigma ^{\ast }(t_{0},x_{0})X\sigma
(t_{0},x_{0})]-\psi
_{i_{0}}(t_{0},x_{0}) && \\
-\left\langle \frac{\gamma }{\varepsilon }(x_{0}-y_{0})|x_{0}-y_{0}|^{2%
\gamma -2}+\frac{\gamma }{\varepsilon }(x_{0}-\overline{x})|x_{0}-\overline{x%
}|^{2\gamma -2}+\theta (2\gamma +2)x_{0}|x_{0}|^{2\gamma
},b(t_{0},x_{0})\right\rangle \leq 0 && \\
\text{ and }d-\frac{1}{2}Tr[\sigma ^{\ast }(t_{0},y_{0})Y\sigma
(t_{0},y_{0})]-(1-\lambda )\psi _{i_{0}}(t_{0},y_{0}) && \\
-\left\langle \frac{\gamma }{\varepsilon }(x_{0}-y_{0})|x_{0}-y_{0}|^{2%
\gamma -2}-\theta (2\gamma +2)y_{0}|y_{0}|^{2\gamma
},b(t_{0},y_{0})\right\rangle \geq 0 &&
\end{eqnarray*}%
which implies that:
\begin{eqnarray}
-c-d &\leq &\frac{1}{2}Tr[\sigma ^{\ast }(t_{0},x_{0})X\sigma
(t_{0},x_{0})-\sigma ^{\ast }(t_{0},y_{0})Y\sigma (t_{0},y_{0})]
\label{viscder} \\
&&+\left\langle \frac{\gamma }{\varepsilon }(x_{0}-y_{0})|x_{0}-y_{0}|^{2%
\gamma -2},b(t_{0},x_{0})-b(t_{0},y_{0})\right\rangle  \notag \\
&&+\left\langle \frac{\gamma }{\varepsilon }(x_{0}-\overline{x})|x_{0}-%
\overline{x}|^{2\gamma -2},b(t_{0},x_{0})\right\rangle +\left\langle
\theta
(2\gamma +2)x_{0}|x_{0}|^{2\gamma },b(t_{0},x_{0})\right\rangle  \notag \\
&&+\left\langle \theta (2\gamma +2)y_{0}|y_{0}|^{2\gamma
},b(t_{0},y_{0})\right\rangle +\psi _{i}(t_{0},x_{0})-(1-\lambda
)\psi _{i}(t_{0},y_{0}).  \notag
\end{eqnarray}%
But from (\ref{derive}) there exist two constants $C$, $C_{1}$ and
$C_{2}$ such that:
\begin{equation*}
||a_{1}(x_{0},y_{0})||\leq C|x_{0}-y_{0}|^{2\gamma -2},\ \
(||a_{2}(x_{0})||\vee ||a_{2}(y_{0})||)\leq C_{1}\theta \,\mbox{ and }%
||a_{3}(x_{0})||\leq \frac{C_{2}}{\varepsilon
}|x_{0}-\overline{x}|^{2\gamma -2}.
\end{equation*}%
As
\begin{equation*}
B=B(t_{0},x_{0},y_{0})=\frac{1}{\varepsilon }%
\begin{pmatrix}
a_{1}(x_{0},y_{0}) & -a_{1}(x_{0},y_{0}) \\
-a_{1}(x_{0},y_{0}) & a_{1}(x_{0},y_{0})%
\end{pmatrix}%
+%
\begin{pmatrix}
a_{2}(x_{0})+a_{3}(x_{0}) & 0 \\
0 & a_{2}(y_{0})%
\end{pmatrix}%
\end{equation*}%
then
\begin{equation*}
B\leq \frac{C}{\varepsilon }|x_{0}-y_{0}|^{2\gamma -2}%
\begin{pmatrix}
I & -I \\
-I & I%
\end{pmatrix}%
+C_{1}\theta I+\frac{C_{2}}{\varepsilon
}|x_{0}-\overline{x}|^{2\gamma -2}I.
\end{equation*}%
It follows that:
\begin{eqnarray}
B+\varepsilon _{1}B^{2} &\leq &C(\frac{1}{\varepsilon }|x_{0}-y_{0}|^{2%
\gamma -2}+\frac{\varepsilon _{1}}{\varepsilon
^{2}}|x_{0}-y_{0}|^{4\gamma
-4})%
\begin{pmatrix}
I & -I \\
-I & I%
\end{pmatrix}
\\
&&+C_{1}\theta I+C_{2}(\frac{1}{\varepsilon
}|x_{0}-\overline{x}|^{2\gamma -2}+\frac{\varepsilon
_{1}}{\varepsilon ^{2}}|x_{0}-\overline{x}|^{4\gamma -4})I  \notag
\end{eqnarray}%
where $C$, $C_{1}$ and $C_{2}$ which hereafter may change from line
to line. Choosing now $\varepsilon _{1}=\varepsilon $, yields the
relation
\begin{eqnarray}
B+\varepsilon _{1}B^{2} &\leq &\frac{C}{\varepsilon
}(|x_{0}-y_{0}|^{2\gamma
-2}+|x_{0}-y_{0}|^{4\gamma -4})%
\begin{pmatrix}
I & -I \\
-I & I%
\end{pmatrix}
\label{ineg_matreciel} \\
&&+C_{1}\theta I+\frac{C_{2}}{\varepsilon
}(|x_{0}-\overline{x}|^{2\gamma -2}+|x_{0}-\overline{x}|^{4\gamma
-4})I.  \notag
\end{eqnarray}%
Now, from (\ref{regbs1}), (\ref{lemmeishii}) and
(\ref{ineg_matreciel}) we get:
\begin{eqnarray*}
&&\frac{1}{2}Tr[\sigma ^{\ast }(t_{0},x_{0})X\sigma
(t_{0},x_{0})-\sigma
^{\ast }(t_{0},y_{0})Y\sigma (t_{0},y_{0})] \\
&\leq &\frac{C}{\varepsilon }(|x_{0}-y_{0}|^{2\gamma
}+|x_{0}-y_{0}|^{4\gamma -2}) \\
&&+C_{1}\theta (1+|x_{0}|^{2}+|y_{0}|^{2})+\frac{C_{2}}{\varepsilon }(|x_{0}-%
\overline{x}|^{2\gamma -2} \\
&&+|x_{0}-\overline{x}|^{4\gamma -4})(1+|x_{0}|^{2}+|y_{0}|^{2}),
\end{eqnarray*}

\begin{equation*}
\left\langle \frac{\gamma }{\varepsilon
}(x_{0}-y_{0})|x_{0}-y_{0}|^{2\gamma
-2},b(t_{0},x_{0})-b(t_{0},y_{0})\right\rangle \leq
\frac{C^{2}}{\varepsilon }|x_{0}-y_{0}|^{2\gamma }.
\end{equation*}%
Next,
\begin{equation*}
\left\langle \frac{\gamma }{\varepsilon }(x_{0}-\overline{x})|x_{0}-%
\overline{x}|^{2\gamma -2},b(t_{0},x_{0})\right\rangle \leq \frac{C}{%
\varepsilon }|x_{0}-\overline{x}|^{2\gamma -3}|x_{0}|
\end{equation*}%
and finally,
\begin{equation*}
\left\langle \theta (2\gamma +2)x_{0}|x_{0}|^{2\gamma
},b(t_{0},x_{0})\right\rangle +\left\langle \theta (2\gamma
+2)y_{0}|y_{0}|^{2\gamma },b(t_{0},y_{0})\right\rangle \leq \theta
C(1+|x_{0}|^{2\gamma +2}+|y_{0}|^{2\gamma +2}).
\end{equation*}%
So that by plugging into (\ref{viscder}) and note that $\lambda >0$
we obtain:
\begin{eqnarray*}
\beta &\leq &\frac{C}{\varepsilon }(|x_{0}-y_{0}|^{2\gamma
}+|x_{0}-y_{0}|^{4\gamma -2})+C_{1}\theta (1+|x_{0}|^{2}+|y_{0}|^{2}) \\
&&+\frac{C_{2}}{\varepsilon }(|x_{0}-\overline{x}|^{2\gamma -2}+|x_{0}-%
\overline{x}|^{4\gamma -4})(1+|x_{0}|^{2}+|y_{0}|^{2})+\frac{C^{2}}{%
\varepsilon }|x_{0}-y_{0}|^{2\gamma } \\
&&+\frac{C}{\varepsilon }|x_{0}-\overline{x}|^{2\gamma
-3}|x_{0}|+\theta
C(1+|x_{0}|^{2\gamma +2}+|y_{0}|^{2\gamma +2})-\frac{\lambda }{t_{0}^{2}}+%
\frac{\gamma }{\varepsilon }(t_{0}-\overline{t})(t_{0}-\overline{t}%
)^{2\gamma -2} \\
&&+\psi _{i_{0}}(t_{0},x_{0})-(1-\lambda )\psi
_{i_{0}}(t_{0},y_{0}).
\end{eqnarray*}%
By sending $\varepsilon \rightarrow 0$, $\lambda \rightarrow 0$,
$\theta \rightarrow 0$ and taking into account of the continuity of
$\psi _{i_{0}}$ and $\gamma \geq 2$, we obtain $\beta \leq 0$ which
is a contradiction. The proof of Theorem \ref{uni} is now complete.
$\Box $

As a by-product we have the following corollary:

\begin{cor}
Let $(v^{1},\ldots ,v^{m})$ be a viscosity solution of
(\ref{sysvi1}) which
satisfies a polynomial growth condition then for $i=1,\ldots ,m$ and $%
(t,x)\in \lbrack 0,T]\times \mathbb{R}^{k}$,
\begin{equation*}
v^{i}(t,x)=\limfunc{esssup}_{(\delta ,\xi )\in
\mathcal{D}_{t}^{i}}E\left[ \int_{t}^{T}\psi
_{u_{s}}(s,X_{s}^{t,x})ds-\sum_{n\geq 1}g_{u_{\tau _{n-1}}u_{\tau
_{n}}}(\tau _{n},X_{\tau _{n}}^{t,x})\mathbf{1}_{[\tau
_{n}<T]}\right] .
\end{equation*}
\end{cor}

\section{Numerical results}

We consider now some numerical examples of the optimal switching problem (%
\ref{sysvi1}).\newline

\begin{exa}
In this example we consider an optimal switching problem with two
modes, where\newline
$T=1$, $b=x$, $\sigma=\sqrt{2}x$, $g_{12}(t,x)=0$, $%
g_{21}(t,x)=0.1|x|+0.5t+2 $, $\psi_1(t,x)=x+0.75t+1$, $\psi_2(t,x)=0.1x+t-1$


\begin{figure}[h]
\begin{center}
\includegraphics{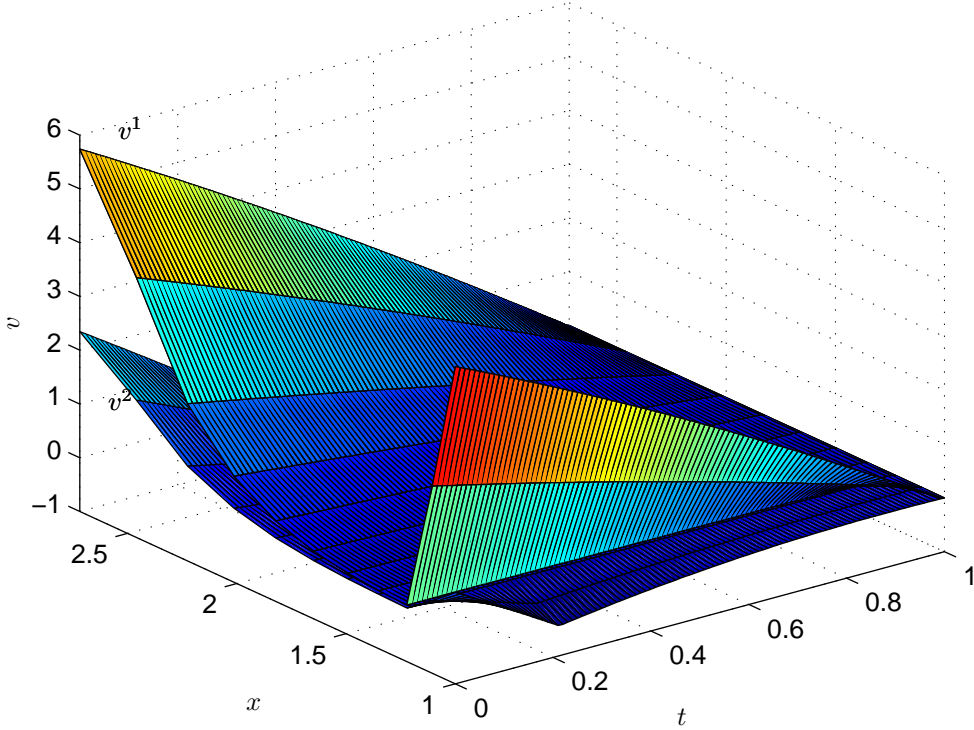}
\end{center}
\caption{Surfaces of $v^1$ and $v^2$.}
\end{figure}
\end{exa}

\begin{exa}
We now consider the case of 3 modes where $T=1$, $b=x$, $\sigma
=\sqrt{2}x$,
$g_{12}(t,x)=0$, $g_{13}(t,x)=0$, $g_{21}(t,x)=|x|+t+4$, $g_{23}(t,x)=0$, $%
g_{31}(t,x)=|x|+t+1$, $g_{32}(t,x)=4t+0.5$, $\psi _{1}(t,x)=x+2t+1$,
$\psi _{2}(t,x)=-x+t-2$ and finally $\psi _{3}(t,x)=-x+t-2$.

\begin{figure}[h]
\begin{center}
\includegraphics{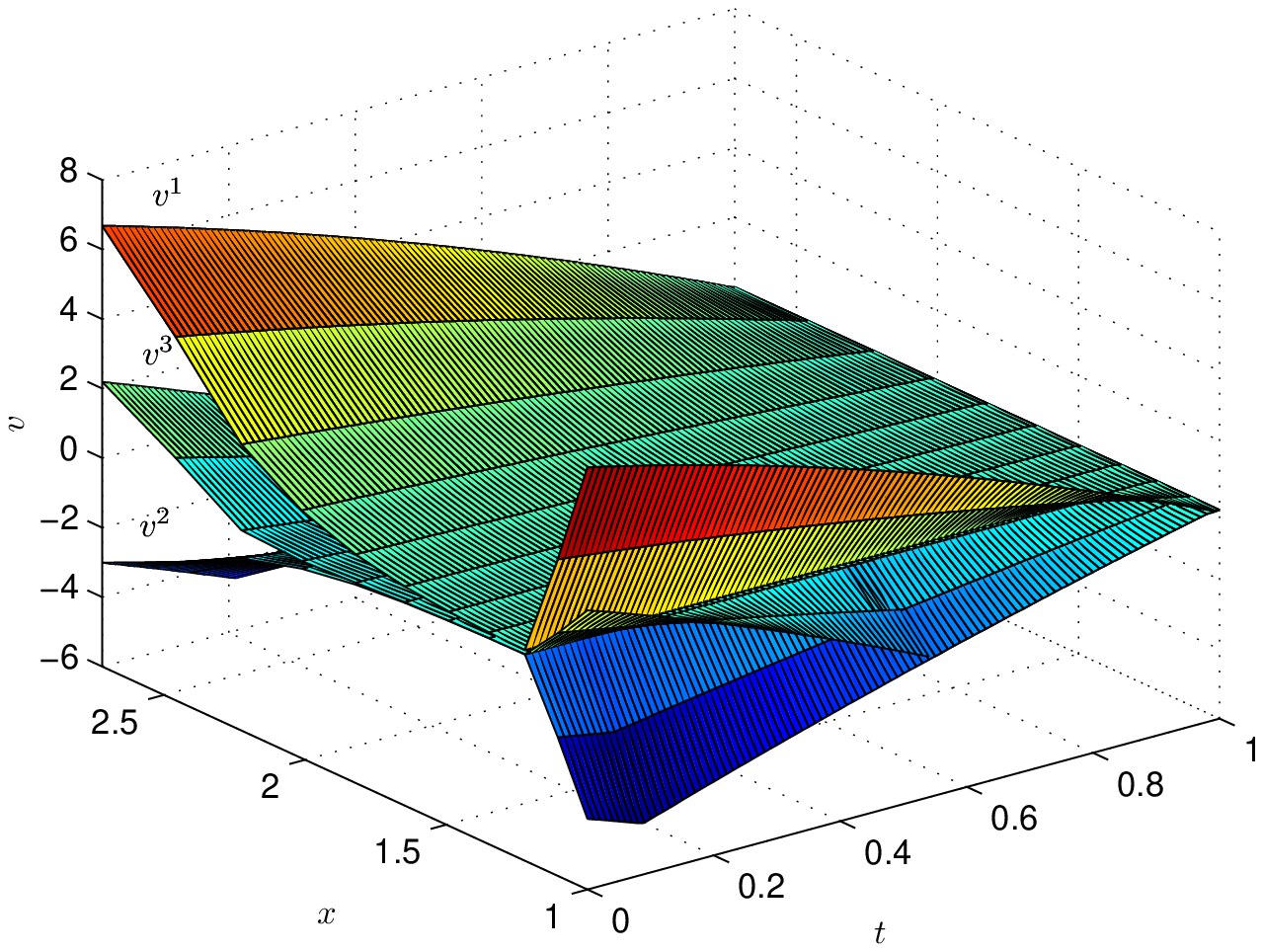}
\end{center}
\caption{Surfaces of $v^{1}$, $v^{3}$ and $v^{2}$.}
\end{figure}
\end{exa}

\section{Appendix}

\indent$Proof$ of Lemma \ref{env-sne}. From (\ref{eqvt}) we have for any $%
i\in \mathcal{I}$ and $0\leq t\leq T$%
\begin{equation}
Y_{t}^{i}=\limfunc{esssup}_{\tau \geq t}E\left[ \int_{t}^{\tau }\psi
_{i}(s,X_{s})ds+\max\limits_{j\in \mathcal{I}^{-i}}(-g_{ij}(\tau
,X_{\tau })+Y_{\tau }^{j})\mathbf{1}_{[\tau
<T]}|\mathcal{F}_{t}\right] .
\end{equation}%
This also means that the process $(Y_{t}^{i}+\int_{0}^{t}\psi
_{i}(s,X_{s})ds)_{0\leq t\leq T}$ is a supermartingale which
dominates
\begin{equation*}
\left( \int_{0}^{t}\psi _{i}(s,X_{s})ds+\max\limits_{j\in \mathcal{I}%
^{-i}}(-g_{ij}(t,X_{t})+Y_{t}^{j})\mathbf{1}_{[t<T]}\right) _{0\leq
t\leq T}.
\end{equation*}%
This implies that the process $(\mathbf{1}_{[u_{\tau _{1}^{\ast
}}=i]}(Y_{t}^{i}+\int_{\tau _{1}^{\ast }}^{t}\psi
_{i}(s,X_{s})ds))_{t\geq \tau _{1}^{\ast }}$ is a supermartingale
which dominates
\begin{equation*}
\left( \mathbf{1}_{[u_{\tau _{1}^{\ast }}=i]}\int_{\tau _{1}^{\ast
}}^{t}\psi _{i}(s,X_{s})ds+\max\limits_{j\in \mathcal{I}%
^{-i}}(-g_{ij}(s,X_{t})+Y_{t}^{j})\mathbf{1}_{[t<T]}\right) _{t\geq
\tau _{1}^{\ast }}.
\end{equation*}%
Since $\mathcal{I}$ is finite, the process $(\sum_{i\in \mathcal{I}}\mathbf{1%
}_{[u_{\tau _{1}^{\ast }}=i]}(Y_{t}^{i}+\int_{\tau _{1}^{\ast
}}^{t}\psi _{i}(s,X_{s})ds))_{t\geq \tau _{1}^{\ast }}$ is also a
supermartingale which dominates
\begin{equation*}
\left( \sum_{i\in \mathcal{I}}\mathbf{1}_{[u_{\tau _{1}^{\ast
}}=i]}\int_{\tau _{1}^{\ast }}^{t}\psi
_{i}(s,X_{s})ds+\max\limits_{j\in
\mathcal{I}^{-i}}(-g_{ij}(t,X_{t})+Y_{t}^{j})\mathbf{1}_{[t<T]}\right)
_{t\geq \tau _{1}^{\ast }}.
\end{equation*}%
Thus, the process $Y_{t}^{u_{\tau _{1}^{\ast }}}+\int_{\tau
_{1}^{\ast }}^{t}\psi _{u_{\tau _{1}^{\ast }}}(s,X_{s})ds)_{t\geq
\tau _{1}^{\ast }}$ is a supermartingale which is greater than
\begin{equation*}
\left( \int_{\tau _{1}^{\ast }}^{t}\psi _{u_{\tau _{1}^{\ast
}}}(s,X_{s})ds+\max\limits_{j\in \mathcal{I}^{-u_{\tau _{1}^{\ast
}}}}(-g_{u_{\tau _{1}^{\ast }}j}(t,X_{t})+Y_{t}^{j})\mathbf{1}%
_{[t<T]}\right) _{t\geq \tau _{1}^{\ast }}.
\end{equation*}%
To complete the proof it remains to show that it is the smallest one
which has this property and use the characterization of the Snell
envelope see e.g. \cite{CK, Elka, ham}.\newline Indeed, let
$(Z_{t})_{0\leq t\leq T}$ be a supermartingale of class $[D]$ such
that, for any $\tau _{1}^{\ast }\leq t\leq T$,
\begin{equation*}
Z_{t}\geq \int_{\tau _{1}^{\ast }}^{t}\psi _{u_{\tau _{1}^{\ast
}}}(s,X_{s})ds+\max\limits_{j\in \mathcal{I}^{-u_{\tau _{1}^{\ast
}}}}(-g_{u_{\tau _{1}^{\ast
}}j}(t,X_{t})+Y_{t}^{j})\mathbf{1}_{[t<T]}.
\end{equation*}%
It follows that for every $\tau _{1}^{\ast }\leq t\leq T$,
\begin{equation*}
\mathbf{1}_{[u_{\tau _{1}^{\ast }}=i]}Z_{t}\geq \mathbf{1}_{[u_{\tau
_{1}^{\ast }}=i]}\left( \int_{\tau _{1}^{\ast }}^{t}\psi
_{i}(s,X_{s})ds+\max\limits_{j\in \mathcal{I}%
^{-i}}(-g_{ij}(t,X_{t})+Y_{t}^{j})\mathbf{1}_{[t<T]}\right) .
\end{equation*}%
But, the process $(\mathbf{1}_{[u_{\tau _{1}^{\ast
}}=i]}Z_{t})_{t\geq \tau _{1}^{\ast }}$ is a supermartingale and for
every $t\geq \tau _{1}^{\ast }$,
\begin{equation*}
\mathbf{1}_{[u_{\tau _{1}^{\ast
}}=i]}Y_{t}^{i}=\limfunc{esssup}_{\tau \geq t}E\left[
\mathbf{1}_{[u_{\tau _{1}^{\ast }}=i]}\left( \int_{t}^{\tau }\psi
_{i}(s,X_{s})ds+\max\limits_{j\in \mathcal{I}^{-i}}(-g_{ij}(\tau
,X_{\tau })+Y_{\tau }^{j})\mathbf{1}_{[\tau <T]}\right)
|\mathcal{F}_{t}\right] .
\end{equation*}%
It follows that, for every $\tau _{1}^{\ast }\leq t\leq T$,
\begin{equation*}
\mathbf{1}_{[u_{\tau _{1}^{\ast }}=i]}Z_{t}\geq \mathbf{1}_{[u_{\tau
_{1}^{\ast }}=i]}(Y_{t}^{i}+\int_{\tau _{1}^{\ast }}^{t}\psi
_{i}(s,X_{s})ds).
\end{equation*}%
Summing over $i$, we get, for every $\tau _{1}^{\ast }\leq t\leq T$,
\begin{equation*}
Z_{t}\geq Y_{t}^{u_{\tau _{1}^{\ast }}}+\int_{\tau _{1}^{\ast
}}^{t}\psi _{u_{\tau _{1}^{\ast }}}(s,X_{s})ds.
\end{equation*}%
Hence, the process $(Y_{t}^{u_{\tau _{1}^{\ast }}}+\int_{\tau
_{1}^{\ast }}^{t}\psi _{u_{\tau _{1}^{\ast }}}(s,X_{s})ds)_{t\geq
\tau _{1}^{\ast }}$ is the Snell envelope of
\begin{equation*}
\left( \int_{\tau _{1}^{\ast }}^{t}\psi _{u_{\tau _{1}^{\ast
}}}(s,X_{s})ds+\max\limits_{j\in \mathcal{I}^{-u_{\tau _{1}^{\ast
}}}}(-g_{u_{\tau _{1}^{\ast }}j}(t,X_{t})+Y_{t}^{j})\mathbf{1}%
_{[t<T]}\right) _{t\geq \tau _{1}^{\ast }},
\end{equation*}%
which completes the proof of the Lemma. $\Box $

\end{document}